\documentclass{article}
\usepackage[english]{babel}
\usepackage{amsmath,amssymb,stmaryrd,enumerate,bbm,latexsym,theorem}
\usepackage[utf8]{inputenc}
\usepackage{graphicx}
\usepackage{float}

\usepackage[top=2.5cm,bottom=2.5cm,right=2cm,left=2cm]{geometry}
\usepackage{authblk}
\numberwithin{equation}{section}

\newcommand{\assign}{:=}
\newcommand{\tmem}[1]{{\em #1\/}}
\newcommand{\tmmathbf}[1]{\ensuremath{\boldsymbol{#1}}}

\newcommand{\nobracket}{}
\newcommand{\tmop}[1]{\ensuremath{\operatorname{#1}}}

\newcommand{\tmtextit}[1]{{\itshape{#1}}}
\newenvironment{itemizedot}{\begin{itemize} }{\end{itemize}}
\newenvironment{proof}{\noindent\textbf{Proof\ }}{\hspace*{\fill}$\Box$\medskip}
\newcounter{nnacknowledgments}

{\theorembodyfont{\rmfamily}\newtheorem{acknowledgments*}[nnacknowledgments]{Acknowledgments}}
\newtheorem{theorem}{Theorem}[section]

\newtheorem{lemma}[theorem]{Lemma}
\newtheorem{proposition}[theorem]{Proposition}
\usepackage[hidelinks]{hyperref}

\date{}
\author{Tej-Eddine Ghoul \thanks{email: teg6@nyu.edu}}
\author{Nader Masmoudi \thanks{email: nm30@nyu.edu}}
\author{Eliot Pacherie  \thanks{email: ep2699@nyu.edu}}
\affil{NYUAD Research Institute, New York University Abu Dhabi}

\begin{document}

\title{Nonlinear enhanced dissipation in viscous Burgers type equations}

\maketitle

\begin{abstract}
  We construct a class of infinite mass functions for which solutions of the
  viscous Burgers equation decay at a better rate than solution of the heat
  equation for the same initial data in this class. In other words, we show an
  enhanced dissipation coming from a nonlinear transport term. We compute the
  asymptotic profile in this class for both equations. For the viscous Burgers
  equation, the main novelty is the construction and description of a time
  dependent profile with a boundary layer, which enhances the dissipation. This profile will be stable
  up to a computable nonlinear correction depending on the perturbation. We
  also extend our results to other convection-diffusion equations.
\end{abstract}

\begin{center}
\end{center}

\section{Introduction and presentation of the results}

We are interested in this paper by the long time behavior of solutions to a
generalized viscous Burgers equation on the real line
\begin{equation}
  \left\{\begin{array}{l}
    \partial_t u - \partial_x^2 u + \partial_x \left( \frac{u^2}{2} + J (u)
    \right) = 0\\
    u_{| t = 0 \nobracket} (x) = \frac{\kappa_{\pm}}{| x |^{\alpha}} (1 + o_{x
    \rightarrow \pm \infty} (1))
  \end{array}\right. \label{gburgers}
\end{equation}
for $\alpha \in] 0, 1 [, \kappa_+, \kappa_- > 0$ and $J$ a smooth function
sastisfying $| J (u) | \leqslant C | u |^3$. Remark that $u_{| t = 0
\nobracket}$ is not integrable, and $J = 0$ correspond to the classical
viscous Burgers equation.

It is well known that for the heat equation $\partial_t u - \partial_x^2 u =
0$, for an initial data $u_0 \in L^1 (\mathbb{R})$ we have the asymptotic
profile
\[ \sqrt{t} u \left( z \sqrt{t}, t \right) \rightarrow \frac{\int_{\mathbb{R}}
   u_0}{\sqrt{4 \pi}} e^{- \frac{z^2}{4}} \]
when $t \rightarrow + \infty$, uniformly in $z \in \mathbb{R}$. A similar
result holds for the viscous Burgers equation $\partial_t u - \partial_x^2 u +
u \partial_x u = 0$ for initial data $u_0 \in L^1 (\mathbb{R})$, (see
{\cite{MR1124296}}, {\cite{MR3583525}}, {\cite{MR2313029}}), as we have
\[ \sqrt{t} u \left( z \sqrt{t}, t \right) \rightarrow \frac{2 \left(
   e^{\frac{M}{2}} - 1 \right) e^{- z^2 / 4}}{e^{\frac{M}{2}} \sqrt{4 \pi} +
   \left( 1 - e^{\frac{M}{2}} \right) \int_{- \infty}^z e^{- s^2 / 4} d s} \]
when $t \rightarrow + \infty$, uniformly in $z \in \mathbb{R}$, where $M =
\int_{\mathbb{R}} u_0$. The same result holds with the term $J (u)$ in the
equation of (\ref{gburgers}).

Although the limit profile is changed by the Burgers term $u \partial_x u$,
the decay rate and the scaling in time are still the same as for the heat
equation. In both case, the $L^{\infty}$ norm of the solution decays like
$t^{- \frac{1}{2}}$.

Other asymptotic behavior results have been established in other
convection-diffusion equations for initial datas in $L^1 (\mathbb{R})$, we
refer to {\cite{arxivzuazuabis}} and references therein, as well as {\cite{MR2666480}},
{\cite{dix_2002},\cite{articleEVZ}},{\cite{MR2369491}}, {\cite{MR2172561}}. See also {\cite{8200fb8cdd914507857ce9e152705d47}} for some results on non integrable initial datas.

\

Going back to our problem (\ref{gburgers}), as a comparison, we look first at
the heat equation for this type of infinite mass initial data. There, up to a
rescaling, we can show that the solution converges to a global attractor.

\begin{proposition}
  \label{heatfo}For $\kappa > 0, \alpha \in] 0, 1 [$, consider $f$ the
  solution of the heat equation $\partial_t f - \partial_x^2 f = 0$ for an
  initial condition $f_0 \in C^0 (\mathbb{R})$ that satisfies
  \[ f_0 (x) = \frac{\kappa (1 + o_{| x | \rightarrow + \infty} (1))}{(1 + | x
     |)^{\alpha}} . \]
  Then, uniformly in $z \in \mathbb{R}$, we have the convergence
  \[ t^{\frac{\alpha}{2}} f \left( \sqrt{t} z, t \right) \rightarrow
     \frac{\kappa}{\sqrt{4 \pi}} \int_{\mathbb{R}} \frac{1}{| y |^{\alpha}}
     e^{- (z - y)^2 / 4} d y \]
  when $t \rightarrow + \infty$.
\end{proposition}

This result is first proven in {\cite{Kamin1985LargeTB}} and the proof is redone in Annex \ref{sec3} to make this paper self contained. Remark that the decay in time of
$f$ is slower than if $f_0$ would be in $L^1 (\mathbb{R})$, in fact, $t^{-
\frac{\alpha}{2}}$ is the size of $t^{- \frac{1}{2}} \int_{-
\sqrt{t}}^{\sqrt{t}} f_0$. Furthermore, the asymptotic profile is smooth, and
behaves like $\kappa | z |^{- \alpha}$ at infinity, connecting back to the
initial data.

\

In this paper, we will construct a stable solution of (\ref{gburgers}) that
converges, up to a rescaling, to an asymptotic profile. However, it will have
two main differences compared to the result of Proposition \ref{heatfo}.
First, the rescaling will not be the same, and surprisingly, the solution will
decay in time like $t^{- \frac{\alpha}{1 + \alpha}}$, that is faster than the
heat equation for the same initial data. The scales of the rescaling are thus dicated by the nonlinear term, which happens also in the nonlinear heat equation with a pure power term, see for instance {\cite{Kamin1985LargeTB}} or {\cite{bookSouplet}}.
Secondly, the asymptotic profile will
be, in the rescaling where it is of size $1$, discontinuous at one point. This discontinuity can be seen as a sort of boundary layer (although there is no boundaries in this problem) that helps the dissipation.
About this enhanced dissipation, we can state the following result.

\begin{theorem}
  \label{airport}Given $\alpha \in \left] \frac{1}{4}, 1 \right[, \kappa_+,
  \kappa_- > 0$, there exists an initial data $u_{| t = 0 \nobracket}$ with
  $u_{| t = 0 \nobracket} (x) = \frac{\kappa_{\pm}}{| x |^{\alpha}} (1 + o_{x
  \rightarrow \pm \infty} (1))$ such that the solution of $\partial_t u -
  \partial_x^2 u + u \partial_x u = 0$ for this initial data satisfies
  \[ t^{\frac{\alpha}{1 + \alpha}} \| u (., t) \|_{L^{\infty} (\mathbb{R})}
     \leqslant c_0 \]
  where $c_0 > 0$ is a constant independent of time. Furthermore, this
  solution is stable in some sense.
\end{theorem}

See Theorem \ref{prop111} for a more precise statement and the shape of the
asymptotic profile, and Proposition \ref{thfinal} for a statement in the case
$J \neq 0$. By Proposition \ref{heatfo}, for the initial data of Theorem
\ref{airport}, if it was instead evolving following the heat equation, we
would have $t^{\frac{\alpha}{1 + \alpha}} \| u (., t) \|_{L^{\infty}
(\mathbb{R})} \rightarrow + \infty$ when $t \rightarrow + \infty$. This means
that the additional nonlinear transport term improved the dissipation.
Enhanced dissipation results are well known for the heat equation with an
additional linear transport term (see for instance {\cite{MR3904158}} {\cite{MR3621820}},
{\cite{MR2434887}}, {\cite{MR4156602}},  {\cite{MR4008523}}, {\cite{MR1225573}} and references
therein) or for Navier-Stokes on $\mathbb{T} \times \mathbb{R}$ (see
{\cite{MR3448924}}, {\cite{MR4134940}}, {\cite{MR4030287}},
{\cite{MR4176913}}). We require $\kappa_+, \kappa_- > 0$ in Theorem
\ref{airport}, and although we can require less, we do not know how to show
that this enhancement is true in general for any $\kappa_+, \kappa_- \in
\mathbb{R}^{\ast}$.

\subsection{Profile for the viscous Burgers equation}\label{ss12}

We focus first in the case $J = 0$ of equation (\ref{gburgers}). There, the
Hopf-Cole formula gives us an explicit formulation of the solution of the
equation. However, since our goal is to be able to generalized it for any $J$,
we will not use it here. The results we can obtain with the Hopf-Cole formula
will be the subject of a companion paper.

Here, we want to construct an approximate solution of the viscous Burgers
equation in the right scaling.

\subsubsection{The rescaled problem and the underlying ODE}\label{plo}

We consider the viscous Bugers equation
\[ \partial_t u - \partial_x^2 u + u \partial_x u = 0 \]
and $\alpha \in] 0, 1 [$. We want to do a change of variable such that the
terms $\partial_t u$ and $u \partial_x u$ are the dominant ones. We define
$\varepsilon (t) \assign t^{\frac{\alpha - 1}{\alpha + 1}}$ and
\[ h (z, \varepsilon (t)) = t^{\frac{\alpha}{1 + \alpha}} u \left( z
   t^{\frac{1}{1 + \alpha}}, t \right), \]
leading to the equivalent equation
\begin{equation}
  \frac{1 - \alpha}{\alpha + 1} \varepsilon \partial_{\varepsilon} h +
  \frac{\alpha}{1 + \alpha} h + \frac{z \partial_z h}{1 + \alpha} - h
  \partial_z h + \varepsilon \partial_z^2 h = 0. \label{key4050}
\end{equation}
Remark that the term coming from the laplacian, $\varepsilon \partial_z^2 h$,
is small when $\varepsilon \rightarrow 0$ (that is $t \rightarrow + \infty$).
This means that at this scale, the nonlinear effect dominates the dynamic.
Interestingly, if we simply remove the laplacian, we get the Burgers equation
$\partial_t u + u \partial_x u = 0$, for which the $L^{\infty}$ norm is
conserved. Since we will show some decay stronger than the heat equation, this
means that although the laplacian is fading out, it still has a major effect
on the dynamic.

\

We want to construct, for $\varepsilon > 0$ small, a solution to the ODE
problem
\begin{equation}
  \left\{\begin{array}{l}
    \frac{\alpha}{1 + \alpha} h + \left( \frac{z}{1 + \alpha} - h \right)
    \partial_z h + \varepsilon \partial_z^2 h = 0\\
    h (z) = \kappa_+ | z |^{- \alpha} (1 + o_{z \rightarrow + \infty} (1))\\
    h (z) = \kappa_- | z |^{- \alpha} (1 + o_{z \rightarrow - \infty} (1))
  \end{array}\right. \label{ODE}
\end{equation}
for $\kappa_+, \kappa_- \in \mathbb{R}^{\ast}$. This will give us an
approximate solution of (\ref{gburgers}). Remark that the problem (\ref{ODE})
is doubly degenerate when $\varepsilon \rightarrow 0$: the coefficient in
front of the term with two derivatives goes to 0, but also, the limit problem
when $\varepsilon = 0$ is ill defined when $h (z) = \frac{z}{1 + \alpha}$,
since then the coefficient in front of the term with one derivative cancels
out.

\subsubsection{The case $\varepsilon = 0$}\label{e0}

We consider in this section the ODE
\[ \left\{\begin{array}{l}
     \frac{\alpha}{1 + \alpha} h + \frac{z \partial_z h}{1 + \alpha} - h
     \partial_z h = 0\\
     h (z_0) = b
   \end{array}\right. \]
for some given $(z_0, b) \in \mathbb{R}^2$ and $\alpha \in] 0, 1 [$. This is
the equation of (\ref{ODE}) with $\varepsilon = 0$. We summarize here the
properties of the solutions.

\begin{figure}[htp]
    \centering
    \includegraphics{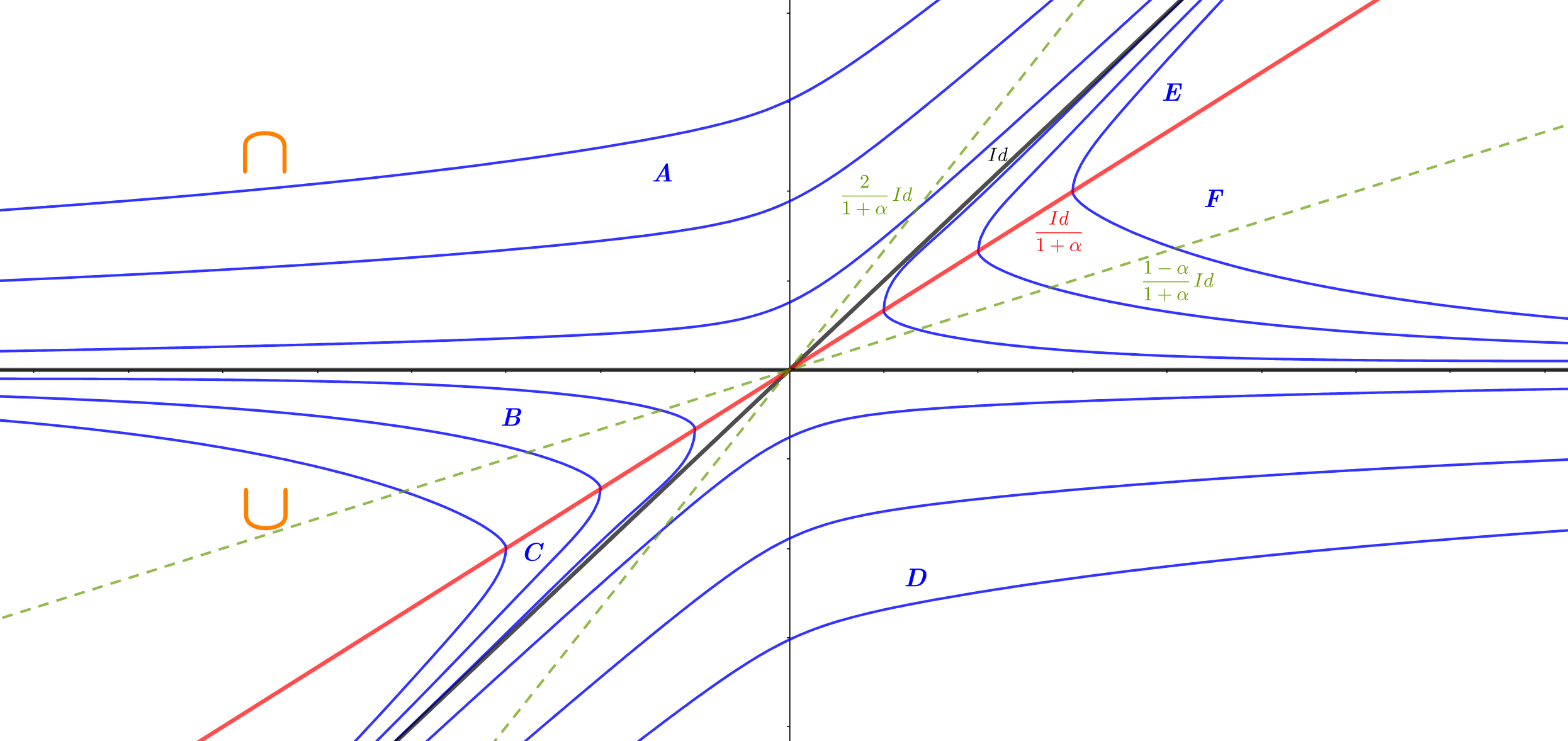}
\end{figure}

\

First, remark that $h (z) = z$ and $h (z) = 0$ are solutions of this equation
(they are the two black lines on the figure). Furthermore, we can write the
equation as
\[ \partial_z h \left( \frac{z}{1 + \alpha} - h \right) + \frac{\alpha}{1 +
   \alpha} h = 0, \]
which is ill defined if $h (z) = \frac{z}{1 + \alpha}$ at some point. This is
the red line in the figure. The blue curves are solution of the equation. In
particular, we have to take $(z_0, b) \in \mathbb{R}^2$ with $b \neq
\frac{z_0}{1 + \alpha}$ for the equation to make sense.

Now, we divide the set
\[ \left\{ (z_0, b) \in \mathbb{R}^2, b \neq \frac{z_0}{1 + \alpha} \right\}
   =\tmmathbf{A} \cup \tmmathbf{B} \cup \mathbf{C} \cup \mathbf{D} \cup
   \mathbf{E} \cup \mathbf{F} \]
with
\[ \tmmathbf{A} \assign \{ (z_0, b) \in \mathbb{R}^2, b > \max (0, z_0) \}, \]
\[ \tmmathbf{B} \assign \left\{ (z_0, b) \in \mathbb{R}^2, 0 > b >
   \frac{z_0}{1 + \alpha} \right\}, \]
\[ \mathbf{C} \assign \left\{ (z_0, b) \in \mathbb{R}^2, \frac{z_0}{1 +
   \alpha} > b > z_0 \right\}, \]
\[ \mathbf{D} \assign \{ (z_0, b) \in \mathbb{R}^2, b < \min (0, z_0) \}, \]
\[ \mathbf{E} \assign \left\{ (z_0, b) \in \mathbb{R}^2, z_0 > b >
   \frac{z_0}{1 + \alpha} \right\}, \]
and
\[ \mathbf{F} \assign \left\{ (z_0, b) \in \mathbb{R}^2, \frac{z_0}{1 +
   \alpha} > b > 0 \right\} . \]

On the figure, the separation between these sets are the red and black lines
(the role of the dotted green lines will be explained later).

The equation has a symmetry: if $z \rightarrow h (z)$ is a solution then so is
$z \rightarrow - h (- z)$. Remark also that since this is a first order ODE,
solutions cannot cross the axes $0, \tmop{Id}$ and $\frac{\tmop{Id}}{1 +
\alpha}$. In particular if a solution has a point in a bold set $\mathbf{J}
\in \{ \tmmathbf{A}, \tmmathbf{B}, \mathbf{C}, \mathbf{D}, \mathbf{E},
\mathbf{F} \}$, then it is fully included in $\mathbf{J}$.

If $(z_0, b) \in \tmmathbf{A}$, then the solution $h$ is global, and
\[ \lim_{z \rightarrow + \infty} h (z) - z = 0, h (z) \sim \kappa | z |^{-
   \alpha} \tmop{when} z \rightarrow - \infty \]
for some $\kappa > 0$ determined by $(z_0, b)$.

If $(z_0, b) \in \tmmathbf{B}$, then the solution is defined on $] - \infty,
z_b [$ for some $z_b > z_0$ determined by $(z_0, b)$, and
\[ \lim_{z \rightarrow z_b} h (z) = \frac{z_b}{1 + \alpha}, \lim_{z
   \rightarrow z_b} h' (z) = - \infty, h (z) \sim - \kappa | z |^{- \alpha}
   \tmop{when} z \rightarrow - \infty \]
for some $\kappa > 0$ determined by $(z_0, b)$.

In both cases, $(z_0, b) \rightarrow \kappa$ is surjective in
$\mathbb{R}_+^{\ast}$.

\

If $(z_0, b) \in \mathbf{C}$, then the solution is defined on $] - \infty,
z_b [$ for some $z_b > z_0$ determined by $(z_0, b)$, and
\[ \lim_{z \rightarrow z_b} h (z) = \frac{z_b}{1 + \alpha}, \lim_{z
   \rightarrow z_b} h' (z) = + \infty, \lim_{z \rightarrow - \infty} h (z) - z
   = 0. \]
By symmetry, we describe similarly the domains and limits if $(z_0, b) \in
\mathbf{D} \cup \mathbf{E} \cup \mathbf{F}$.

\

In particular, remark that they are no continuous solutions to the problem
\[ \left\{\begin{array}{l}
     \frac{\alpha}{1 + \alpha} h + \frac{z \partial_z h}{1 + \alpha} - h
     \partial_z h = 0\\
     h (z) = \kappa_+ | z |^{- \alpha} (1 + o_{z \rightarrow + \infty} (1))\\
     h (z) = \kappa_- | z |^{- \alpha} (1 + o_{z \rightarrow - \infty} (1))
   \end{array}\right. \]
for $\kappa_+, \kappa_- \in \mathbb{R}^{\ast}$. Therefore, we expect solution
of (\ref{ODE}) to have jumps in the limit $\varepsilon \rightarrow 0$.

In the next subsection, we will give some conditions to describe what jumps
are possible in the limit $\varepsilon \rightarrow 0$. We will show that for
$\kappa_+, \kappa_- > 0$, at most one is a viscosity solution.

\subsubsection{Viscosity solutions}\label{ss1235}

First, if $h$ is a solution of (\ref{ODE}) with $\varepsilon = 0$ \ in the
distribution sense, then it must satisfy the Rankine-Hugoniot conditions.
Here, it states that at any discontinuity $z_c \in \mathbb{R}$, we must have
\begin{equation}
  \frac{1}{2} (h (z_c^+) + h (z_c^-)) = \frac{z_c}{1 + \alpha} .
\end{equation}
On the figure, this means that the middle point of any jump must be on the red
line. This prevents for instance jumps one bold set to itself, but also for
instance from $\mathbf{F}$ to $\mathbf{D}$.

The dotted green line $\frac{2}{1 + \alpha} \tmop{Id}$ and $\frac{1 -
\alpha}{1 + \alpha} \tmop{Id}$ are the ones such that the red line is the
middle between $0$ and $\frac{2}{1 + \alpha} \tmop{Id}$, and the middle
between $\tmop{Id}$ and $\frac{1 - \alpha}{1 + \alpha} \tmop{Id}$. Remark that
for any $\alpha \in] 0, 1 [$, we have the order $0 < \frac{1 - \alpha}{1 +
\alpha} < \frac{1}{1 + \alpha} < 1 < \frac{2}{1 + \alpha}$.

We continue, remark that if $\frac{\alpha}{1 + \alpha} h_{\varepsilon} +
\frac{z \partial_z h_{\varepsilon}}{1 + \alpha} - h_{\varepsilon} \partial_z
h_{\varepsilon} + \varepsilon \partial_z^2 h_{\varepsilon} = 0$ and
$\partial_z h_{\varepsilon} (z) = 0$ for some $z \in \mathbb{R}$, then
\[ \partial_z^2 h_{\varepsilon} (z) = \frac{- \alpha}{(1 + \alpha)
   \varepsilon} h_{\varepsilon} (z) . \]
This is represented by the two orange cups on the figure: if $h_{\varepsilon}
(z) > 0, h_{\varepsilon}' (z) = 0$, then $h''_{\varepsilon} (z) < 0$, so near
$z$ the function $h$ looks like the cup.

This means that, if we expect $h$, a solution of (\ref{ODE}) with $\varepsilon
= 0$ with some discontinuities, to be the limit when $\varepsilon \rightarrow
0$ of a sequence of function $h_{\varepsilon}$ solution of (\ref{ODE}), since
the $h_{\varepsilon}$ are smooth, then some jumps cannot happen. For instance,
although the Rankine-Hugoniot conditions allows jumps from $\mathbf{F}$ to
$\mathbf{E}$, they are not viscous (this would require the existence of $z$
such that $h_{\varepsilon} (z) > 0, h_{\varepsilon}' (z) = 0$ and
$h''_{\varepsilon} (z) \geqslant 0$).

\

We continue. We infer that it is not possible to have two jumps that cross
the axis $\{ z = 0 \}$. Indeed, otherwise we denote $z_a < z_b$ two
consecutives values such that $h_{\varepsilon} (z_a) = h_{\varepsilon} (z_b) =
0$, and integrating the equation between $z_a$ and $z_b$ leads to
\[ \frac{\alpha - 1}{\alpha + 1} \int_{z_a}^{z_b} h_{\varepsilon} (s) d s +
   \varepsilon (h_{\varepsilon}' (z_b) - h'_{\varepsilon} (z_a)) = 0, \]
but this is impossible since either $h_{\varepsilon}' (z_b) - h'_{\varepsilon}
(z_a) > 0$ and $h_{\varepsilon} > 0$ on $[z_a, z_b]$, or $h_{\varepsilon}'
(z_b) - h'_{\varepsilon} (z_a) < 0$ and $h_{\varepsilon} < 0$ on $[z_a, z_b]$.

Finally, suppose that a solution has a point $z_a < 0$ where $h (z_a) =
\frac{2}{1 + \alpha} z_a, h' (z_a) > 0$ and a point $z_b > z_a$ such that $h
(z_b) = 0$. Then, integrating the equation between $z_a$ and $z_b$ leads to
\[ \frac{\alpha - 1}{\alpha + 1} \int_{z_a}^{z_b} h_{\varepsilon} (s) d s +
   \varepsilon (h_{\varepsilon}' (z_b) - h'_{\varepsilon} (z_a)) = 0. \]
When $\varepsilon \rightarrow 0$, this also lead to a contradiction. This
prevent the possibility to have a solutions having jumps between
$\tmmathbf{B}$ and $\mathbf{D}$ followed by a jump from $\mathbf{D}$ to
$\tmmathbf{A}$.

We summarize these conditions in the following image.

\begin{figure}[H]
    \centering
    \includegraphics[width=20cm]{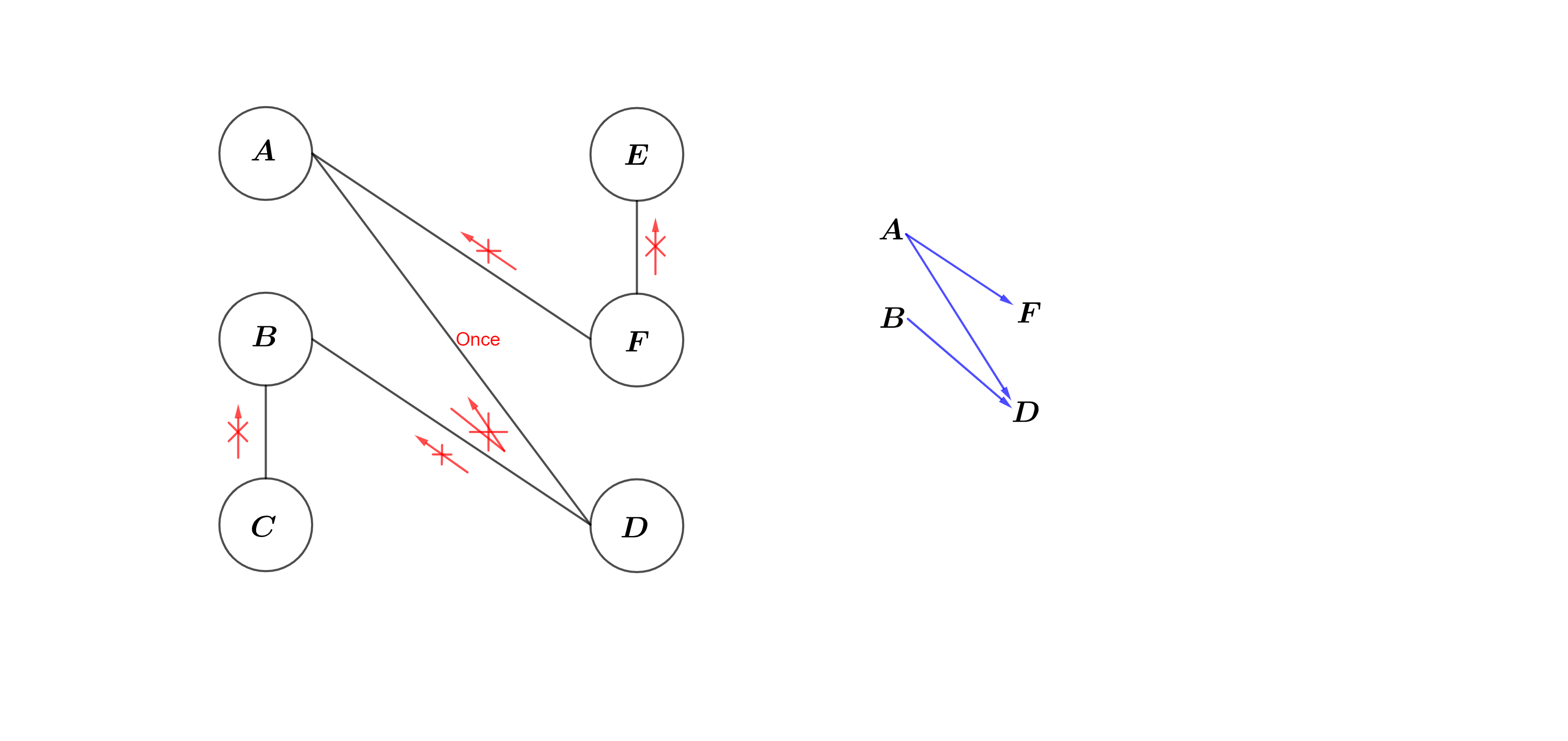}
\end{figure}

Jumps are not possible from a bold set to itself. Two bolds sets are connected
by a black line if there is a possible jump between them satisfying the
Rankine-Hugoniot condition. Crossed red arrows are jumps that are forbidden by
the viscosity conditions. The jump between $\tmmathbf{A}$ and $\mathbf{D}$ can
only be done once.

\

These viscosity conditions severely limit what jumps are allowed. We are
looking for solutions starting from $\tmmathbf{A}$ or $\tmmathbf{B}$ and
ending in $\mathbf{D}$ or $\mathbf{F}$. These conditions imposes that, for the
cases $\tmmathbf{A} \rightarrow \mathbf{F}, \tmmathbf{A} \rightarrow
\mathbf{D}$ and $\tmmathbf{B} \rightarrow \mathbf{D}$, only a single jump is
possible. Furthermore, in these cases (that is $\kappa_+, \kappa_-$ having the
same sign or $\kappa_- > 0, \kappa_+ < 0$), the position of the jump is fully
determined by $\kappa_+$ and $\kappa_-$. For instance if $\kappa_+, \kappa_- >
0$, these two values gives on which blue curves in $\tmmathbf{A}$ and
$\mathbf{F}$ the solution is, and we can check that there is only one value
$z_c$ for which the red line is the middle point of these two blue curves.

\

For the last case $\tmmathbf{B} \rightarrow \mathbf{F}$, where it seems that
no connection is possible, we omitted the case where the jump does not finish
in a bold set, but finish on the identity line. It is therefore maybe possible
to go from $\tmmathbf{B}$ to $\mathbf{F}$ with two jumps, both connecting to
$\{ (z, z), z \in \mathbb{R} \}$. It is however difficult to prove or disprove
that such a thing might happen. It might also be possible that the solution
of the viscous Burgers equation with such an initial data simply does not
converge with this rescaling.

\

The results of these three above subsections are shown in subsections
\ref{lane8} to \ref{ifonly}.

\subsubsection{Construction of the profile for small $\varepsilon > 0$}

In this section, given $\kappa_+, \kappa_- > 0, \alpha \in] 0, 1 [$ and
$\varepsilon > 0$ small enough, we want to construct a solution of the ODE
problem
\begin{equation}
  \left\{\begin{array}{l}
    \frac{\alpha}{1 + \alpha} \mathbbmss{h}_{\varepsilon} + \frac{z \partial_z
    \mathbbmss{h}_{\varepsilon}}{1 + \alpha} -\mathbbmss{h}_{\varepsilon} \partial_z
    \mathbbmss{h}_{\varepsilon} + \varepsilon \partial_z^2
    \mathbbmss{h}_{\varepsilon} = 0\\
    \mathbbmss{h}_{\varepsilon} (z) = \kappa_{\pm} | z |^{- \alpha} (1 + o_{z
    \rightarrow \pm \infty} (1)) .
  \end{array}\right. \label{14pl}
\end{equation}
But first, we define the function $h_0$ by the unique solution to the problem
\begin{equation}
  \left\{\begin{array}{l}
    \frac{\alpha}{1 + \alpha} h_0 + \frac{z \partial_z h_0}{1 + \alpha} - h_0
    \partial_z h_0 = 0\\
    h_0 (z) = \kappa_+ | z |^{- \alpha} (1 + o_{z \rightarrow + \infty} (1))
  \end{array}\right. \label{renegade1}
\end{equation}
if $z > z_c$, and the unique viscous solution to
\begin{equation}
  \left\{\begin{array}{l}
    \frac{\alpha}{1 + \alpha} h_0 + \frac{z \partial_z h_0}{1 + \alpha} - h_0
    \partial_z h_0 = 0\\
    h_0 (z) = \kappa_- | z |^{- \alpha} (1 + o_{z \rightarrow - \infty} (1))
  \end{array}\right. \label{renegade2}
\end{equation}
if $z < z_c$, where $z_c \in \mathbb{R}$ is the position of the jump given by
the conditions described above, uniquely determined by $\kappa_+, \kappa_-,
\alpha$. We will show in subsection \ref{ifonly} the existence and uniqueness
of the solution of these problems.

We define
\[ h_0 (z_c^{\pm}) \assign \lim_{\nu \rightarrow 0} h_0 (z_c \pm \nu), \]
as the function $h_0$ is discontinous at $z_c$.

\

We construct a solution of (\ref{14pl}) using a shooting method, and this
solution will be close to $h_0$ far from $z_c$, see the following result.

\begin{proposition}
  \label{profileODEd}For any $\kappa > 0, \alpha \in] 0, 1 [$, there exists
  $z_c \in \mathbb{R}, \varepsilon_0 > 0$ such that, for $\varepsilon_0 >
  \varepsilon > 0$, there exists two $C^1$ functions $\varepsilon \rightarrow
  z_c (\varepsilon), a (\varepsilon)$ with
  \[ z_c (\varepsilon) \rightarrow z_c, a (\varepsilon) \rightarrow \frac{h_0
     (z_c^-) - h_0 (z_c^+)}{2} \]
  when $\varepsilon \rightarrow 0$, such that the solution of the ODE problem
  \[ \left\{\begin{array}{l}
       \frac{\alpha}{1 + \alpha} \mathbbmss{h}_{\varepsilon} + \left( \frac{z}{1
       + \alpha} -\mathbbmss{h}_{\varepsilon} \right) \partial_z
       \mathbbmss{h}_{\varepsilon} + \varepsilon \partial_z^2
       \mathbbmss{h}_{\varepsilon} = 0\\
       \mathbbmss{h}_{\varepsilon} (z_c (\varepsilon)) = \frac{z_c
       (\varepsilon)}{1 + \alpha}, \mathbbmss{h}_{\varepsilon}' (z_c
       (\varepsilon)) = \frac{- a (\varepsilon)^2}{2 \varepsilon}
     \end{array}\right. \]
  satisfies
  \[ \mathbbmss{h}_{\varepsilon} (z) = \kappa_{\pm} | z |^{- \alpha} (1 + o_{z
     \rightarrow \pm \infty} (1)) . \]
  Furthermore, there exists $w_0 > 0$ depending only on $\alpha$ and $\kappa$,
  such that
  \[ \left\| \mathbbmss{h}_{\varepsilon} (z) - h_0 (z_c^+) - \frac{2 a
     (\varepsilon) e^{- a (\varepsilon) \left( \frac{z - z_c
     (\varepsilon)}{\varepsilon} \right)}}{1 + e^{- a (\varepsilon) \left(
     \frac{z - z_c (\varepsilon)}{\varepsilon} \right)}} \right\|_{L^{\infty}
     \left( \left[ z_c (\varepsilon) - w_0 \varepsilon \ln
     \frac{1}{\varepsilon}, z_c (\varepsilon) + w_0 \varepsilon \ln
     \frac{1}{\varepsilon} \right] \right)} \rightarrow 0 \]
  and
  \[ \| (1 + | z |)^{\alpha} (\mathbbmss{h}_{\varepsilon} - h_0) (z)
     \|_{L^{\infty} \left( \mathbb{R} \backslash \left[ z_c (\varepsilon) -
     w_0 \varepsilon \ln \frac{1}{\varepsilon}, z_c (\varepsilon) + w_0
     \varepsilon \ln \frac{1}{\varepsilon} \right] \right)} \rightarrow 0 \]
  when $\varepsilon \rightarrow 0$. Finally, we have $| \partial_{\varepsilon}
  z_c (\varepsilon) | + | \partial_{\varepsilon} a (\varepsilon) | \leqslant K
  \left( \ln \frac{1}{\varepsilon} \right)^2$.
\end{proposition}

Section \ref{Daisy} is devoted to the proof of this result. For $\varepsilon
\neq 0$ and $a = \frac{h_0 (z_c^-) - h_0 (z_c^+)}{2}$, the solution to
\[ \left\{\begin{array}{l}
     \frac{\alpha}{1 + \alpha} h_{\varepsilon} + \left( \frac{z}{1 + \alpha} -
     h_{\varepsilon} \right) \partial_z h_{\varepsilon} + \varepsilon
     \partial_z^2 h_{\varepsilon} = 0\\
     h_{\varepsilon} (z_c) = \frac{z_c}{1 + \alpha}, h_{\varepsilon}' (z_c
     (\varepsilon)) = \frac{- a^2}{2 \varepsilon}
   \end{array}\right. \]
does not satisfies $h_{\varepsilon} (z) = \kappa_{\pm} | z |^{- \alpha} (1 +
o_{z \rightarrow \pm \infty} (1))$, but
\[ h_{\varepsilon} (z) = (\kappa_{\pm} + o_{\varepsilon \rightarrow 0} (1)) |
   z |^{- \alpha} (1 + o_{z \rightarrow \pm \infty} (1)) . \]
This is why, to get the exact same equivalent at $+ \infty$, we need to change
slightly $z_c (\varepsilon)$ and $a (\varepsilon)$. We use the notation
$\mathbbmss{h}_{\varepsilon}$ for solutions of the problem (\ref{14pl}) (that is
depending on the behavior at $\pm \infty$) and $h_{\varepsilon}$ for solutions
depending on its value at $z_c$.

\

The function $\mathbbmss{h}_{\varepsilon}$ in Proposition \ref{profileODEd} is
close to the discontinuous function $h_0$ except in a vicinity of $z_c$, the
discontinuity point. $\mathbbmss{h}_{\varepsilon}$ solves (\ref{ODE}) but not
(\ref{key4050}) because of the term $\frac{1 - \alpha}{1 + \alpha} \varepsilon
\partial_{\varepsilon} \mathbbmss{h}_{\varepsilon}$, however this term is small
compared to the other ones. We will show the stability of
$\mathbbmss{h}_{\varepsilon}$ in a space that contains in particular this error
term in the next subsection.

\

To construct $\mathbbmss{h}_{\varepsilon}$, we found the right scale around
$z_c$ to have now a continuous function (it is $\frac{z - z_c
(\varepsilon)}{\varepsilon} \simeq 1$ rather than $z \simeq 1$). The proof of
Proposition \ref{profileODEd} is done in two parts. First, we compute the
first order in $\varepsilon$ of the solution in $\left[ z_c (\varepsilon) -
w_0 \varepsilon \ln \frac{1}{\varepsilon}, z_c (\varepsilon) + w_0 \varepsilon
\ln \frac{1}{\varepsilon} \right]$ for some $w_0 > 0$ large but independent of
$\varepsilon$, and we show that at the boundaries of this interval, it become
close to the value of $h_0$ at the same point. Then, outside this interval,
$h_0$ and $\mathbbmss{h}_{\varepsilon}$ verify similar equation for small
$\varepsilon$, and start with similar values. We thus show that they stay
close.

\subsubsection{Stability of the profile}

We recall that $\varepsilon (t) = t^{\frac{\alpha - 1}{\alpha + 1}}$ and
$\mathbbmss{h}_{\varepsilon}$ is the solution described in Proposition
\ref{profileODEd}. We want to show that if at a time $T > 0$ large, we solve
the viscous Burgers equation with the initial data $\mathbbmss{h}_{\varepsilon
(T)} + f_0$ at time $T$ in the rescaled variables, then for all times $t
\geqslant T$ we stay close to $\mathbbmss{h}_{\varepsilon (t)}$. Interestingly,
$\mathbbmss{h}_{\varepsilon (t)}$ will not be the first order, we need to modify
it nonlinearly, depending on $f_0$. It turns out that the mass of $f_0$ will
change the profile near $z_c$, in a non negligeable way. The stability result
is as follow.

\begin{theorem}
  \label{prop111}Given $\alpha \in \left] \frac{1}{4}, 1 \right[, \kappa_+,
  \kappa_- > 0$, there exists $T_0 > 0$ such that, for any $T \geqslant T_0$,
  there exists $\nu > 0$ depending on $T$ such that, considering
  $\mathbbmss{h}_{\varepsilon (t)}, z_c (t), a (t)$ defined in Proposition
  \ref{profileODEd}, the solution $u$ to the problem
  \[ \left\{\begin{array}{l}
       \partial_t u - \partial_x^2 u + u \partial_x u = 0\\
       u_{| t = T \nobracket} (x) = T^{- \frac{\alpha}{1 + \alpha}}
       \mathbbmss{h}_{\varepsilon (T)} \left( T^{- \frac{1}{1 + \alpha}} x
       \right) + f_0 (x)
     \end{array}\right. \]
  with $f_0 \in H^2 (\mathbb{R})$ and
  \[ \| (1 + | x |)^3 f_0 (x) \|_{L^{\infty} (\mathbb{R})} + \| \partial_x f_0
     \|_{H^1 (\mathbb{R})} + \left| \int_{\mathbb{R}} f_0 \right| \leqslant
     \nu \]
  satisfies for any $t \geqslant T$ that
  \[ \left\| t^{\frac{\alpha}{1 + \alpha}} u (x, t) -\mathbbmss{h}_{\varepsilon
     (t)} \left( t^{- \frac{1}{1 + \alpha}} x \right) -\mathfrak{u} \left(
     t^{\frac{1 - \alpha}{1 + \alpha}} \left( x t^{- \frac{1}{1 + \alpha}} -
     z_c (t) \right) \right) \right\|_{L^{\infty} (\mathbb{R})} = o_{t
     \rightarrow + \infty} (1), \]
  where $\mathfrak{u}$ is the unique solution to the problem
  \[ \left\{\begin{array}{l}
       - \partial_x \mathfrak{u}+ \left( \frac{a (t) (e^{- a (t) x} - 1)}{1 +
       e^{- a (t) x}} \right) \mathfrak{u}+ \frac{\mathfrak{u}^2}{2} = 0\\
       \int_{\mathbb{R}} \mathfrak{u}_{} = \int_{\mathbb{R}} f_0 .
     \end{array}\right. \]
\end{theorem}

Section \ref{Amywiles} is devoted to the proof of this result. Let us make
some remarks about it.
\begin{itemizedot}
  \item This result implies Theorem \ref{airport} and gives us the asymptotic
  profile
  \[ t^{\frac{\alpha}{1 + \alpha}} u \left( z t^{\frac{1}{1 + \alpha}}, t
     \right) \rightarrow h_0 (z) \]
  when $t \rightarrow + \infty$ for any $z \neq z_c$. The convergence is
  uniform on $\mathbb{R}$ if we remove any open set containing $z_c$. In a
  vicinity of $z_c$ we still have convergence to some limit, and there this
  limit depends on $\mathfrak{u}_{}$, that is $f_0$. Remark that
  $\mathfrak{u}$ depends nonliearly on $f_0$, and thus the correction coming
  from $\mathfrak{u}$ is not simply a modulation on the parameters of
  $\mathbbmss{h}_{\varepsilon (t)}$ (that is $\varepsilon (t), a (t)$ or $z_c
  (t)$), even if for small values of $\int_{\mathbb{R}} f_0$, we have $\mathfrak{u}
  \simeq \partial_{z_c} \mathbbmss{h}_{\varepsilon (t)}$.
  
  \item With the conditions on $f_0$, we check that our initial data
  \[ T^{- \frac{\alpha}{1 + \alpha}} \mathbbmss{h}_{\varepsilon (T)} \left( T^{-
     \frac{1}{1 + \alpha}} x \right) + f_0 (x) \]
  decays like $\kappa_{\pm} | x |^{- \alpha}$ when $x \rightarrow \pm \infty$,
  and $f_0$ is small when compared to the main profile, since $\nu$ depends on
  $T$. Also, the condition $\alpha > \frac{1}{4}$ is a technical one, we
  expect the result to hold for any $\alpha \in] 0, 1 [$. This condition will
  be used to show that $\partial_{\varepsilon} \mathbbmss{h}_{\varepsilon}$ has
  enough decay at $\pm \infty$ to estimate it in $H^1 (\mathbb{R})$, see
  subsection \ref{magitek}. 
\end{itemizedot}

The core idea of the proof is to write the solution for $t \geqslant T$ as
\[ u (x, t) = t^{- \frac{\alpha}{1 + \alpha}} \left( \mathbbmss{h}_{\varepsilon
   (t)} \left( t^{- \frac{1}{1 + \alpha}} x \right) +\mathfrak{u} \left(
   t^{\frac{1 - \alpha}{1 + \alpha}} \left( x t^{- \frac{1}{1 + \alpha}} - z_c
   (t) \right) \right) + f (x, t) \right), \]
and now the error $f$ is massless (that is $\int_{\mathbb{R}} f = 0$). We
write it $f = \partial_x g$, and it turns out that we can integrate the
equation to have a new equation on $g$. We show there some coercivity on the
linear part on $g$ in $H^2 (\mathbb{R})$, and we control the nonlinear part,
from which we deduce that $\| g \|_{H^2 (\mathbb{R})} \rightarrow 0$ when $t
\rightarrow + \infty$.

\subsection{Generalisation to equation (\ref{gburgers})}\label{ss13}

Our approach also work for the equation $\partial_t u - \partial_x^2 u +
\partial_x \left( \frac{u^2}{2} + J (u) \right) = 0$, if $J$ satisfies
\[ | J (x) | + | x J' (x) | + | x^2 J'' (x) | \leqslant C_0 | x |^3 \]
for some $C_0 > 0$.

\begin{proposition}
  \label{thfinal}For $\alpha \in \left] \frac{2}{3}, 1 \right[$ and $T_0, \nu$
  depending on $C_0$, the result of Theorem \ref{prop111} also holds for the
  problem
  \[ \left\{\begin{array}{l}
       \partial_t u - \partial_x^2 u + \partial_x \left( \frac{u^2}{2} + J (u)
       \right) = 0\\
       u_{| t = T \nobracket} (x) = T^{- \frac{\alpha}{1 + \alpha}}
       \mathbbmss{h}_{\varepsilon (T)} \left( T^{- \frac{1}{1 + \alpha}} x
       \right) + f_0 (x) .
     \end{array}\right. \]
\end{proposition}

Section \ref{gxd} is devoted to the proof of this result. It is done simply by
checking that the term $\partial_x (J (u))$ can be considered as en error term
in the proof of the stability of Theorem \ref{prop111}. This is true because
at the scale where we see the profile $\mathbbmss{h}_{\varepsilon}$, this term is
small compared to the other ones. As before, the condition $\alpha >
\frac{2}{3}$ is a technical one, and is here to make sure that $J (u)$ has
enough decay at $\pm \infty$ to estimate it in $H^1 (\mathbb{R})$.

\subsection{Some related open problems}

Our results should extend easily for values of $\kappa_+, \kappa_- \in
\mathbb{R}^{\ast}$ except in the case $\kappa_+ < 0, \kappa_+ > 0$. There, it
is maybe possible to construct a specific solution, but it is likely that it
is an unstable one. If we generalized to the equation $\partial_t u -
\partial_x^2 u + u^k \partial_x u = 0$ for $k \in \mathbb{N}^{\ast}$, it is
likely that a similar result can be shown with some improvements in the
proofs.

It seems for now difficult to generalize this result to higher dimension, but
it would be of interest, in particular if similar profiles can be constructed
for the 2d Euler or Navier-Stokes equation.

\begin{acknowledgments*}
  The authors are supported by Tamkeen under the NYU Abu Dhabi Research
  Institute grant CG002. The authors have no competing interests to declare that are relevant to the content of this article.
\end{acknowledgments*}

\section{Construction of the profile $\mathbbmss{h}_{\varepsilon}$}\label{Daisy}

This section is devoted to the proof of Proposition \ref{profileODEd} and the
viscosity properties described in the introduction. First, in section
\ref{lane8} we set the change of scaling and compute the Rankine-Hugoniot
condition for the viscous Burgers equation. Section \ref{lane9} is devoted to
the case $\varepsilon = 0$. Section \ref{ifonly} is about the construction of
$h_0$ (which will be the limit of $\mathbbmss{h}_{\varepsilon}$ when $\varepsilon
\rightarrow 0$), as well as the study of its properties. Section
\ref{Timezone} and \ref{Timezone2} are the study of the shooting problem at
the heart of Proposition \ref{profileODEd}, respectively close and far from
the shooting point $z_c$. Section \ref{Timezone3} regroups all these elements
and concludes the proof of Proposition \ref{profileODEd}.

\subsection{Change of variable and viscosity conditions}\label{lane8}

In this subsection, our goal is to prove some results of subsections \ref{plo}
to \ref{ss1235} of the introduction.

\subsubsection{Computation of the underlying ODE problem}

We consider here the equation
\[ \partial_t u - \partial_x^2 u + u \partial_x u = 0. \]
We define
\[ g (z, t) = t^{\frac{\alpha}{1 + \alpha}} u \left( z t^{\frac{1}{1 +
   \alpha}}, t \right), \]
we have
\[ t \partial_t g = \frac{\alpha}{1 + \alpha} g + \frac{z}{1 + \alpha}
   t^{\frac{\alpha}{1 + \alpha} + \frac{1}{1 + \alpha}} \partial_x u +
   t^{\frac{\alpha}{1 + \alpha} + 1} \partial_t u, \]
and
\[ \partial_z g = t \partial_x u, \partial_z^2 g = t^{1 + \frac{1}{1 +
   \alpha}} \partial_x^2 u. \]
Therefore,
\[ t \partial_t g = \frac{\alpha}{1 + \alpha} g + \frac{z}{1 + \alpha}
   t^{\frac{\alpha}{1 + \alpha} + \frac{1}{1 + \alpha} - 1} \partial_z g +
   t^{\frac{\alpha}{1 + \alpha} + 1} \left( t^{- 1 - \frac{1}{1 + \alpha}}
   \partial_z^2 g - t^{- 1 - \frac{\alpha}{1 + \alpha}} g \partial_z g
   \right), \]
that is
\[ t \partial_t g = \frac{\alpha}{1 + \alpha} g + \frac{z \partial_z g}{1 +
   \alpha} - g \partial_z g + t^{\frac{\alpha - 1}{1 + \alpha}} \partial_z^2
   g. \]
We define $\varepsilon (t) = t^{\frac{\alpha - 1}{\alpha + 1}}$ and we do the
change of variable
\[ h (z, \varepsilon) = g (z, t), \]
since
\[ t \partial_t \varepsilon = \frac{\alpha - 1}{\alpha + 1} \varepsilon
   \tmop{we} \tmop{have} \]
\begin{equation}
  \frac{1 - \alpha}{1 + \alpha} \varepsilon \partial_{\varepsilon} h +
  \frac{\alpha}{1 + \alpha} h + \frac{z \partial_z h}{1 + \alpha} - h
  \partial_z h + \varepsilon \partial_z^2 h = 0. \label{toda}
\end{equation}
By Proposition \ref{heatfo}, this scaling is not adapted to the heat equation
with the same initial condition. Remark if we tried to use this scale anyway,
we would get the same equation (\ref{toda}) but without the term $- h
\partial_z h$. When $\varepsilon \rightarrow 0$, the limit problem will be
$\frac{\alpha}{1 + \alpha} h + \frac{z \partial_z h}{1 + \alpha} = 0$, which
has only the solution $h = C z^{- \alpha}$ for some $C > 0$, which is
unbounded at $z = 0$.

\subsubsection{Rankine-Hugoniot condition}

For the equation $\frac{\alpha}{1 + \alpha} h + \frac{z \partial_z h}{1 +
\alpha} - h \partial_z h = 0$, integrating it between $z_c - \nu$ and $z_c +
\nu$ leads, after some computations, to
\begin{eqnarray*}
  &  & \frac{\alpha - 1}{1 + \alpha} \int_{z_c - \nu}^{z_c + \nu} h +
  \frac{1}{1 + \alpha} ((z_c + \nu) h (z_c + \nu) - (z_c - \nu) h (z_c -
  \nu))\\
  & - & \frac{1}{2} (h^2 (z_c + \nu) - h^2 (z_c - \nu))\\
  & = & 0.
\end{eqnarray*}
Therefore, letting $\nu \rightarrow 0$ leads to
\[ \frac{z_c}{1 + \alpha} (h (z_c^+) - h (z_c^-)) - \frac{1}{2} (h^2 (z_c^+) -
   h^2 (z_c^-)) = 0, \]
that we factorize as
\[ (h (z_c^+) - h (z_c^-)) \left( \frac{z_c}{1 + \alpha} - \frac{1}{2} (h
   (z_c^+) + h (z_c^-)) \right) = 0, \]
which is the Rankine-Hugoniot condition stated in the introduction.

\subsection{Some properties of solutions of $\frac{\alpha}{1 + \alpha} h +
\left( \frac{\tmop{Id}}{1 + \alpha} - h \right) h' = 0$}\label{lane9}

Take $(z_0, b) \in \mathbb{R}^2$ with $b \neq \frac{z_0}{1 + \alpha}$, and we
consider here the problem
\begin{equation}
  \left\{\begin{array}{l}
    \frac{\alpha}{1 + \alpha} h + \left( \frac{z}{1 + \alpha} - h \right)
    \partial_z h = 0\\
    h (z_0) = b,
  \end{array}\right. \label{bb1}
\end{equation}
that is the problem described in section \ref{e0}. The fact that if $(z_0, b)
\in \mathbf{J} \in \{ \tmmathbf{A}, \tmmathbf{B}, \mathbf{C}, \mathbf{D},
\mathbf{E}, \mathbf{F} \}$ implies that $(z, h (z)) \in \mathbf{J}$ for all
values of $z$ on which the solution $h$ of (\ref{bb1}) is well defined is a
consequence of standard Cauchy theory arguments, since the boundary between
two bold set is either a solution of (\ref{bb1}), or the set $\left\{ (z_0, b)
\in \mathbb{R}^2, b = \frac{z_0}{1 + \alpha} \right\}$, on which $\partial_z
h$ explode.

\subsubsection{The case $(z_0, b) \in \tmmathbf{A}$}

\begin{lemma}
  \label{16BL}The solution $h$ of (\ref{bb1}) with $(z_0, b) \in \tmmathbf{A}$
  is defined on $\mathbb{R}$ and satisfies
  \[ \lim_{z \rightarrow + \infty} h (z) - z = 0. \]
\end{lemma}

We leave the study of the behavior when $z \rightarrow - \infty$ for
subsection \ref{ifonly}.

\begin{proof}
  We consider $(z_0, b) \in \tmmathbf{A}= \{ (z_0, b) \in \mathbb{R}^2, b >
  \max (0, z_0) \}$. As long as the solution $h$ of (\ref{bb1}) for this
  initial condition exists, we have $(z, h (z)) \in \tmmathbf{A}$, therefore
  $h (z) > \max (0, z)$. We denote $] z_-, z_+ [$ the maximum domain of
  existence of $h$ with $z_-, z_+ \in \mathbb{R} \cup \{ \pm \infty \}$ (by
  definition we have $z_0 \in] z_-, z_+ [$).
  
  \
  
  Suppose that $z_+ \neq + \infty$. There exists $C_0 > 0$ depending on $z_0$
  and $b$ such that, for $z \in [z_0, z_+ [$ we have $\left| \frac{z}{1 +
  \alpha} - h (z) \right| \geqslant C_0 (1 + | z |)$. Indeed, we have $h (z) >
  \max (0, z)$ and $\frac{1}{1 + \alpha} < 1$. In particular, $\frac{z}{1 +
  \alpha} - h (z) \neq 0$ on $[z_0, z_+ [$ and since $h (z) > 0$ we have
  \[ \left| \frac{\partial_z h}{h} \right| \leqslant \frac{K}{(1 + | z |)} \]
  on $[z_0, z_+ [$. With $z_+ < + \infty$, we deduce that $h$ and $\partial_z
  h$ are bounded near $z_+$, which is a contradiction, therefore $z_+ = +
  \infty$. We define for $z \geqslant z_0$ the function $u (z) = h (z) - z >
  0$. It satisfies the equation
  \begin{equation}
    \left( \frac{\alpha}{1 + \alpha} - \partial_z h (z) \right) u + \left(
    \frac{z}{1 + \alpha} - h (z) \right) \partial_z u = 0 \label{b2402}
  \end{equation}
  on $[z_0, + \infty [$. Now, we compute, using the equation satified by $h$,
  that
  \[ \left( \frac{\alpha}{1 + \alpha} - \partial_z h (z) \right) \left(
     \frac{z}{1 + \alpha} - h (z) \right) = \frac{\alpha}{1 + \alpha} \left(
     \frac{z}{1 + \alpha} - h (z) \right) + \frac{\alpha}{1 + \alpha} h (z) =
     \frac{\alpha z}{(1 + \alpha)^2}, \]
  hence
  \[ \frac{\alpha}{1 + \alpha} - \partial_z h (z) = \frac{\alpha z}{(1 +
     \alpha)^2 \left( \frac{z}{1 + \alpha} - h (z) \right)}, \]
  and we can write equation (\ref{b2402}) as
  \[ \frac{\alpha z}{(1 + \alpha)^2} u + \left( \frac{z}{1 + \alpha} - h (z)
     \right)^2 \partial_z u = 0. \]
  Using $\left| \frac{z}{1 + \alpha} - h (z) \right| \geqslant C_0 (1 + | z
  |)$, we deduce that for $z \geqslant \max (1, z_0)$ we have
  \[ \frac{\partial_z u}{u} \leqslant \frac{- C_1}{z} \]
  for some $C_1 > 0$, therefore $u (z) \leqslant K z^{- C_1}$ for $z \geqslant
  \max (1, z_0)$ for some constant $K > 0$, hence $u (z) \rightarrow 0$ when
  $z \rightarrow + \infty$, leading to $h (z) - z \rightarrow 0$ when $z
  \rightarrow + \infty$.
  
  \
  
  On $] z_-, z_0]$, by similar arguments as previously, we have $h > 0$ and
  $\frac{z}{1 + \alpha} - h (z) < 0$. Therefore, on $] z_-, z_0]$ we have
  $\partial_z h (z) > 0, h (z) > 0$, hence $z_- = - \infty$.
\end{proof}

\subsubsection{The case $(z_0, b) \in \mathbf{F}$}

\begin{lemma}
  The solution $h$ of (\ref{bb1}) with $(z_0, b) \in \mathbf{F}$ is defined on
  $] z_-, + \infty [$ for some $z_- > 0$.
\end{lemma}

We also leave the study of the behavior when $z \rightarrow + \infty$ for
subsection \ref{ifonly}.

\begin{proof}
  We consider here $(z_0, b) \in \mathbf{F}= \left\{ (z_0, b) \in
  \mathbb{R}^2, \frac{z_0}{1 + \alpha} > b > 0 \right\}$ in the problem
  (\ref{bb1}). As in the previous subsection, we consider the largest interval
  on which the solution is defined, that we write $] z_-, z_+ [$. We have $z_0
  \in] z_-, z_+ [$ and for $z \in] z_-, z_+ [$, we have
  \[ \frac{z}{1 + \alpha} > h (z) > 0. \]
  A consequence of this and equation $\frac{\alpha}{1 + \alpha} h + \left(
  \frac{z}{1 + \alpha} - h \right) \partial_z h = 0$ is that $\partial_z h <
  0$ on $] z_-, z_+ [$, and with $h > 0$, we deduce that $z_+ = + \infty$. We
  also see that $z_- > 0$ because the condition $\frac{z}{1 + \alpha} > h (z)
  > 0$ can no longer hold at $z = 0$.
\end{proof}

\subsubsection{The remaining cases}

For $(z_0, b) \in \mathbf{E}$, we can deal with the limit for large $z$ as in
the case of $\tmmathbf{A}$, and we can show that $z_- > 0$ as in the case of
$\mathbf{F}$. By symmetry, we show similar properties in $\tmmathbf{B},
\mathbf{C}$ and $\mathbf{D}$.

\subsection{Definition and properties of the profile $h_0$}\label{ifonly}

The goal of this subsection is to show that, given $\alpha \in] 0, 1 [$ and
$\kappa_+, \kappa_- > 0$, there exists a unique value of $z_c$ an a unique
viscous solution of (\ref{renegade1})-(\ref{renegade2}) in the sense of the
introduction. We will also study its properties.

\subsubsection{A connected implicit problem}

We look for an implicit solution of $\frac{\alpha}{1 + \alpha} h + \left(
\frac{z}{1 + \alpha} - h \right) \partial_z h = 0$ of the form $z = g (h)$.
Differentiating with respect to $z$, we have $1 = \partial_z h g' (h)$ and
replacing, we deduce that
\[ g' (h) = \frac{- g (h)}{\alpha h} + \frac{1 + \alpha}{\alpha} . \]
The solution of this equation is of the form $g (h) = h + \frac{\kappa}{| h
|^{\alpha}}$.

\

This is why we define for $\alpha \in] 0, 1 [$ and $\kappa > 0$ of the
function
\[ g_{\kappa} (y) \assign y + \frac{\kappa}{| y |^{\alpha}} . \]
We are interested in the solutions of the implicit problem $z = g_{\kappa} (y
(z))$. First, remark that $g_{\kappa} (y) \rightarrow + \infty$ when $y
\rightarrow 0^{\pm}, g_{\kappa} (y) \rightarrow \pm \infty$ when $y
\rightarrow \pm \infty$ and $g_{\kappa}' (y) \rightarrow 1$ when $y
\rightarrow \pm \infty$. We compute for $y \neq 0$ that
\[ g_{\kappa}' (y) = 1 - \frac{\kappa \alpha}{y | y |^{\alpha}} . \]
In particular, $g_{\kappa}' > 0$ on $] - \infty, 0 [$. We have $g_{\kappa}'
(y) = 0$ if and only if $y = y_{\kappa} = (\kappa \alpha)^{\frac{1}{1 +
\alpha}} .$ This implies that on $[y_0, + \infty [,$we have $g_{\kappa}' (y) >
0$. We compute easily that
\[ g_{\kappa} (y_{\kappa}) = \kappa^{\frac{1}{1 + \alpha}} \left(
   \alpha^{\frac{1}{1 + \alpha}} + \alpha^{- \frac{\alpha}{1 + \alpha}}
   \right) > 0. \]
By the implicit function theorem, given $\kappa_+, \kappa_- > 0$ we construct
two particular branches of functions. First, a smooth function $y^{\ast}_- :
\mathbb{R} \rightarrow] - \infty, 0 [$ solution of $z = g_{\kappa_-}
(y^{\ast}_- (z))$ for any $z \in \mathbb{R}$, defined as the inverse of the
invertible function $g_{\kappa_-} :] - \infty, 0 [\rightarrow \mathbb{R}$, and
another smooth function
\[ y^{\ast}_+ :] g_{\kappa_+} (y_{\kappa_+}), + \infty [\rightarrow]
   y_{\kappa_+}, + \infty [ \]
solution of $z = g_{\kappa_+} (y^{\ast}_+ (z))$, defined as the inverse of
\[ g_{\kappa_+} :] y_{\kappa_+}, + \infty [\rightarrow] g_{\kappa_+}
   (y_{\kappa_+}), + \infty [. \]
We define here $h_{\pm} (z) \assign \frac{\kappa_{\pm}}{| y_{\pm}^{\ast} (z)
|^{\alpha}} = z - y_{\pm}^{\ast} (z)$. Since $g_{\kappa}' (y) \rightarrow 1$
when $| y | \rightarrow + \infty$, we have that $y_{\pm}^{\ast} (z)
\rightarrow z$ when $z \rightarrow \pm \infty$ and therefore $h_{\pm} (z) \sim
\frac{\kappa_{\pm}}{| z |^{\alpha}}$ when $| z | \rightarrow + \infty$. Let us
show that these functions are solution of $\frac{\alpha}{1 + \alpha} h +
\left( \frac{z}{1 + \alpha} - h \right) \partial_z h = 0$.

\begin{lemma}
  \label{L41n}The functions $h_{\pm}$ satisfy, on their domain of definitions,
  the equation
  \[ \frac{\alpha}{1 + \alpha} h_{\pm} + \left( \frac{z}{1 + \alpha} - h_{\pm}
     \right) \partial_z h_{\pm} = 0. \]
\end{lemma}

\begin{proof}
  We first check that $g_{\kappa}' (y) = 1 - \frac{\kappa \alpha}{y | y
  |^{\alpha}} = (1 + \alpha) y - \alpha g_{\kappa} (y)$, and since
  $g_{\kappa_{\pm}} (y^{\ast}_{\pm} (z)) = z$, we have
  \[ y^{\ast}_{\pm} (z) g'_{\kappa_{\pm}} (y^{\ast}_{\pm} (z)) = (1 + \alpha)
     y^{\ast}_{\pm} (z) - \alpha z. \]
  Furthermore, differentiating with respect to $z$ the equation $z =
  g_{\kappa} (y (z))$ leads to
  \[ \partial_z y^{\ast}_{\pm} (z) g_{\kappa_{\pm}}' (y^{\ast}_{\pm} (z)) = 1,
  \]
  therefore
  \[ \partial_z y^{\ast}_{\pm} (z) \left( y^{\ast}_{\pm} (z) - \frac{\alpha}{1
     + \alpha} z \right) = \frac{y^{\ast}_{\pm} (z)}{1 + \alpha} . \]
  Since $h_{\pm} (z) = z - y_{\pm}^{\ast} (z)$, we have
  \[ \left( \frac{z}{1 + \alpha} - h_{\pm} \right) \partial_z h_{\pm} = \left(
     \frac{z}{1 + \alpha} - h_{\pm} \right) (1 - \partial_z y_{\pm}^{\ast}) =
     \frac{z}{1 + \alpha} - h_{\pm} - \partial_z y_{\pm}^{\ast} \left(
     y^{\ast}_{\pm} (z) - \frac{\alpha}{1 + \alpha} z \right) \]
  and thus
  \[ \left( \frac{z}{1 + \alpha} - h_{\pm} \right) \partial_z h_{\pm} =
     \frac{z}{1 + \alpha} - h_{\pm} - \frac{y^{\ast}_{\pm} (z)}{1 + \alpha} =
     \frac{- \alpha}{1 + \alpha} h_{\pm} . \]
\end{proof}

This results allows us to complete the study of the problem (\ref{bb1}) for
$(z_0, b) \in \tmmathbf{A}$ and $(z_0, b) \in \tmmathbf{F}$.

\begin{lemma}
  \label{SvD}Given $(z_0, b) \in \tmmathbf{A}$, the solution of the problem
  (\ref{bb1}) is $h_-$ defined above for the value $\kappa_- = b | b - z_0
  |^{\alpha}$. It satisfies in particular
  \[ \lim_{z \rightarrow - \infty} | z |^{\alpha} h_- (z) = b | b - z_0
     |^{\alpha} . \]
  Furthermore, for $z \leqslant z_0, n \in \mathbb{N}$, there exists $C_n > 0$
  depending on $n, z_0, \alpha$ such that
  \[ | \partial_z^n h_- (z) | \leqslant \frac{C_n}{(1 + | z |)^{\alpha + n}} .
  \]
  Similarly, for $(z_0, b) \in \tmmathbf{F}$, the solution of the problem
  (\ref{bb1}) is $h_+$ defined above for the value $\kappa_+ = b | b - z_0
  |^{\alpha} .$ It satisfies
  \[ \lim_{z \rightarrow + \infty} | z |^{\alpha} h_+ (z) = b | b - z_0
     |^{\alpha} \]
  and for $z \geqslant z_0$,
  \[ | \partial_z^n h_+ (z) | \leqslant \frac{C_n}{(1 + | z |)^{\alpha + n}} .
  \]
\end{lemma}

\begin{proof}
  Given $(z_0, b) \in \tmmathbf{A}$, we look for a value $\kappa$ such that
  \[ z_0 = g_{\kappa} (z_0 - b) \]
  with $g_{\kappa} (z) = z + \frac{\kappa}{| z |^{\alpha}}$. We check that
  then $\kappa = b | b - z_0 |^{\alpha}$. Defining $h_-$ as above, remark that
  $h_- (z_0) = b$ and $h_-$ satisfies $\frac{\alpha}{1 + \alpha} h_- + \left(
  \frac{z}{1 + \alpha} - h_- \right) \partial_z h_- = 0.$ It is therefore the
  solution of (\ref{bb1}) for this choice of $(z_0, b) \in \tmmathbf{A}$. The
  proof is identical if $(z_0, b) \in \tmmathbf{F}$ with $\kappa_+ = b | b -
  z_0 |^{\alpha}$, to the difference that we have to check that $z_0 - b =
  y^{\ast} (z_0) \in] y_{\kappa_+}, + \infty [$, that is
  \[ z_0 - b \geqslant (\kappa \alpha)^{\frac{1}{1 + \alpha}} = | z_0 - b
     |^{\frac{\alpha}{1 + \alpha}} b^{\frac{1}{1 + \alpha}} \alpha^{\frac{1}{1
     + \alpha}} . \]
  This inequality is a consequence of the fact that $(z_0, b) \in
  \tmmathbf{F}$, which implies that $\frac{z_0}{1 + \alpha} > b > 0$, and thus
  $z_0 - b \geqslant \alpha b$.
  
  \
  
  Now, concerning the computations of $\partial_z^n h_{\pm}$, we have
  \[ \partial_z h_{\pm} = \frac{\frac{- \alpha}{1 + \alpha}
     h_{\pm}}{\frac{z}{1 + \alpha} - h_{\pm}}, \]
  hence
  \[ | \partial_z h_{\pm} (z) | \leqslant \frac{C_1}{(1 + | z |)^{\alpha + 1}}
  \]
  and we can conclude by induction.
\end{proof}

We complete this subsection with a technical lemma on the dependency on $b$
and $z_0$ of $h_{\pm}$.

\begin{lemma}
  \label{Vegas}The function $(b, z_0) \rightarrow h_{\pm}$ is differentiable,
  and there exists $K > 0$ depending on $b, z_0$ such that
  \[ \left| h_{\pm} (z) - \frac{b | b - z_0 |^{\alpha}}{| z |^{\alpha}}
     \right| \leqslant \frac{K}{(1 + | z |)^{1 + 2 \alpha}} \]
  as well as
  \[ \left| \partial_b h_{\pm} (z) - \frac{\partial_b (b | b - z_0
     |^{\alpha})}{| z |^{\alpha}} \right| \leqslant \frac{K}{(1 + | z |)^{1 +
     2 \alpha}} \]
  on the domain of definition of $h_{\pm}$.
\end{lemma}

\begin{proof}
  Take $y (z)$ a function solution of the implicit problem
  \[ y (z) + \frac{\kappa}{| y (z) |^{\alpha}} = z \]
  defined on some interval $[z_0, + \infty [$. By the remarks above Lemma
  \ref{L41n} we have $y (z) \sim z$ when $z \rightarrow + \infty$. Writing $y
  = z + \bar{y} (z)$, we check that
  \[ | z |^{\alpha} \bar{y} (z) \left| 1 + \frac{\bar{y} (z)}{z}
     \right|^{\alpha} = - \kappa, \]
  hence
  \[ \bar{y} (z) \sim \frac{- \kappa}{| z |^{\alpha}} \]
  when $z \rightarrow + \infty$. We deduce that $h (z) = \frac{\kappa}{| y (z)
  |^{\alpha}}$ satisfies
  \[ h (z) - \frac{\kappa}{| z |^{\alpha}} \sim \frac{\alpha \kappa^2}{| z
     |^{1 + 2 \alpha}} \]
  when $z \rightarrow + \infty$. We then check easily that we have similarly
  \[ \partial_{\kappa} h (z) - \frac{1}{| z |^{\alpha}} \sim \frac{2 \alpha
     \kappa}{| z |^{1 + 2 \alpha}} \]
  when $z \rightarrow + \infty$, which implies the result of the lemma for
  $h_+$, and a similar proof works for $h_-$.
\end{proof}

\subsubsection{Connection between the jump and the limits at $\pm \infty$}

We recall the notation
\[ h (z_c^{\pm}) = \lim_{\nu \rightarrow 0} h (z_c \pm \nu) . \]

\begin{lemma}
  \label{L42}Take $\alpha \in] 0, 1 [,$ $z_c > 0$ and $\frac{z_c}{1 + \alpha}
  > a > \frac{\alpha}{1 + \alpha} z_c$. Then, there exists a unique function
  $h \in C^{\infty} (\mathbb{R} \backslash \{ z_c \},] 0, + \infty [)$
  solution to the problem
  \[ \left\{\begin{array}{l}
       \frac{\alpha}{1 + \alpha} h + \left( \frac{z}{1 + \alpha} - h \right)
       \partial_z h = 0\\
       h (z_c^+) = \frac{z_c}{1 + \alpha} - a\\
       h (z_c^-) = \frac{z_c}{1 + \alpha} + a.
     \end{array}\right. \]
  Furthermore, there exists $\kappa_{\pm} (z_c, a) > 0$ such that
  \[ h (z) = \frac{\kappa_{\pm} (z_c, a)}{| z |^{\alpha}} (1 + o_{z
     \rightarrow \pm \infty} (1)), \]
  and
  \[ (z_c, a) \rightarrow (\kappa_+ (z_c, a), \kappa_- (z_c, a)) \]
  is a smooth function and a bijection from $\left\{ (z_c, a) \in] 0, + \infty
  [^2, \frac{z_c}{1 + \alpha} > a > \frac{\alpha}{1 + \alpha} z_c \right\}$ to
  $] 0, + \infty [^2$.
\end{lemma}

\begin{proof}
  By Lemma \ref{16BL}, the solution of
  \[ \left\{\begin{array}{l}
       \frac{\alpha}{1 + \alpha} h + \left( \frac{z}{1 + \alpha} - h \right)
       \partial_z h = 0\\
       h (z_c) = \frac{z_c}{1 + \alpha} - a
     \end{array}\right. \]
  with $z_c > 0$ and $0 < h (z_c) < \frac{z_c}{1 + \alpha}$ (that is
  $\frac{z_c}{1 + \alpha} > a > 0$) is well defined for all $z \geqslant z_c$,
  and we have $g_{\kappa_+ (z_c, a)} (z - h (z)) = z$ for
  \[ \kappa_+ (z_c, a) \assign h (z_c) | z_c - h (z_c) |^{\alpha} = \left(
     \frac{z_c}{1 + \alpha} - a \right) \left| \frac{\alpha}{1 + \alpha} z_c +
     a \right|^{\alpha} . \]
  Similarly, the solution of
  \[ \left\{\begin{array}{l}
       \frac{\alpha}{1 + \alpha} h + \left( \frac{z}{1 + \alpha} - h \right)
       \partial_z h = 0\\
       h (z_c) = \frac{z_c}{1 + \alpha} + a
     \end{array}\right. \]
  with $h (z_c) > z_c$ (that is $a > \frac{\alpha}{1 + \alpha} z_c$) is well
  defined for all $z \leqslant z_c$, and we have $g_{\kappa_- (z_c, a)} (z - h
  (z)) = z$ for
  \[ \kappa_- (z_c, a) \assign h (z_c) | z_c - h (z_c) |^{\alpha} = \left(
     \frac{z_c}{1 + \alpha} + a \right) \left| \frac{\alpha}{1 + \alpha} z_c -
     a \right|^{\alpha} . \]
  We deduce that
  \[ (z_c, a) \rightarrow (\kappa_+ (z_c, a), \kappa_- (z_c, a)) \]
  is a smooth function from $\left\{ (z_c, a) \in \mathbb{R}, z_c > 0, a \in
  \left] \frac{\alpha}{1 + \alpha} z_c, \frac{1}{1 + \alpha} z_c \right[
  \right\}$ to $] 0, + \infty [^2$. Let us show that it is a bijection.
  Writing $a = z_c b$ with $b \in \left] \frac{\alpha}{1 + \alpha}, \frac{1}{1
  + \alpha} \right[$, we have
  \[ \kappa_+ = z_c^{1 + \alpha} \left( \frac{1}{1 + \alpha} - b \right)
     \left| \frac{\alpha}{1 + \alpha} + b \right|^{\alpha} \]
  and
  \[ \kappa_- = z_c^{1 + \alpha} \left( \frac{1}{1 + \alpha} + b \right)
     \left| b - \frac{\alpha}{1 + \alpha} \right|^{\alpha}, \]
  hence
  \[ \frac{\kappa_+}{\kappa_-} = \zeta (b) \assign \frac{\left( \frac{1}{1 +
     \alpha} - b \right) \left| b + \frac{\alpha}{1 + \alpha}
     \right|^{\alpha}}{\left( \frac{1}{1 + \alpha} + b \right) \left| b -
     \frac{\alpha}{1 + \alpha} \right|^{\alpha}} . \]
  $\zeta$ is a smooth function of $b$ in $\left] \frac{\alpha}{1 + \alpha},
  \frac{1}{1 + \alpha} \right[$, and
  \[ \lim_{b \rightarrow \left( \frac{\alpha}{1 + \alpha} \right)^+} \zeta (b)
     = + \infty, \lim_{b \rightarrow \left( \frac{1}{1 + \alpha} \right)^-}
     \zeta (b) = 0. \]
  We check that
  \[ \zeta' (b) = \frac{- 2 b \left| b + \frac{\alpha}{1 + \alpha}
     \right|^{\alpha - 1} \left| b - \frac{\alpha}{1 + \alpha} \right|^{\alpha
     - 1}}{\left( \left( \frac{1}{1 + \alpha} + b \right) \left| b -
     \frac{\alpha}{1 + \alpha} \right|^{\alpha} \right)^2} < 0 \]
  hence $\zeta$ is a bijection from $\left] \frac{\alpha}{1 + \alpha},
  \frac{1}{1 + \alpha} \right[$ to $\mathbb{R}^{+ \ast}$. This completes the
  proof of the lemma.
\end{proof}

We now can construct the function $h_0$: for $\kappa_+, \kappa_- > 0$, take
$(z_c (\kappa_+, \kappa_-), a (\kappa_+, \kappa_-)) \in] 0, + \infty [^2$ such
that $\kappa_+ (z_c (\kappa_+, \kappa_-), a (\kappa_+, \kappa_-)) = \kappa_+,
\kappa_- (z_c (\kappa_+, \kappa_-), a (\kappa_+, \kappa_-)) = \kappa_-$, then
$h_0$ is the solution of Lemma \ref{L42} for these values. It is almost a
solution of (\ref{ODE}), but it is discontinuous at $z_c$. It satisfies the
Rankine-Hugoniot condition, and by Lemma \ref{L42} it is the only solution
among the ones behaving like $\kappa_{\pm} | z |^{- \alpha}$ at $\pm \infty$
doing so.

\

Our goal in the next subsection is to construct a better approximation
$\mathbbmss{h}_{\varepsilon}$, that will be continuous at $z_c$, and be close to
$h_0$ away from it when $\varepsilon$ is small.

\subsection{Shooting from $z_c$ and shape of the profile near
it}\label{Timezone}

We want to understand the solution to the problem
\begin{equation}
  \left\{\begin{array}{l}
    \frac{\alpha}{1 + \alpha} h_{\varepsilon} + \left( \frac{z}{1 + \alpha} -
    h_{\varepsilon} \right) \partial_z h_{\varepsilon} + \varepsilon
    \partial_z^2 h_{\varepsilon} = 0\\
    h_{\varepsilon} (z_c) = \frac{z_c}{1 + \alpha}, h'_{\varepsilon} (z_c) =
    \frac{- a^2}{2 \varepsilon}
  \end{array}\right. \label{41pb}
\end{equation}
for some given parameters $z_c, a > 0$ with $\frac{z_c}{1 + \alpha} > a >
\frac{\alpha}{1 + \alpha} z_c$ that for now are independent of $\varepsilon$.
In this subsection, we take $h$ to be the solution of Lemma \ref{L42}
associated to the values of $z_c$ and $a$. The function $h$ is discontinuous
at $z_c$. We want to show that for a right choice of $z_c$ and $a$,
$h_{\varepsilon}$ is close to $h$ far from $z_c$, and we want to compute the
shape of $h_{\varepsilon}$ near $z_c$.

\subsubsection{Estimates in $\left] z_c, z_c + w_0 \varepsilon \ln
\frac{1}{\varepsilon} \right]$}

\begin{lemma}
  \label{cigaretteandkerosine}There exists $w_0 > 0$ depending only on
  $\alpha, z_c, a$ such that $h_{\varepsilon}$, the solution of (\ref{41pb}),
  is well defined on $\left] z_c, z_c + w_0 \varepsilon \ln
  \frac{1}{\varepsilon} \right]$ and satisfies
  \[ \left\| h_{\varepsilon} (z) - h (z) - \frac{2 a e^{- a \left( \frac{z -
     z_c}{\varepsilon} \right)}}{1 + e^{- a \left( \frac{z - z_c}{\varepsilon}
     \right)}} \right\|_{L^{\infty} \left( \left] z_c, z_c + w_0 \varepsilon
     \ln \frac{1}{\varepsilon} \right] \right)} \rightarrow 0 \]
  as well as
  \[ \left\| h'_{\varepsilon} (z) + \frac{2 a^2 e^{- a \left( \frac{z -
     z_c}{\varepsilon} \right)}}{\varepsilon \left( 1 + e^{- a \left( \frac{z
     - z_c}{\varepsilon} \right)} \right)^2} \right\|_{L^{\infty} \left(
     \left] z_c, z_c + w_0 \varepsilon \ln \frac{1}{\varepsilon} \right]
     \right)} \leqslant C \]
  for some constant $C > 0$ depending only on $\alpha, z_c, a$ when
  $\varepsilon \rightarrow 0$, and also
  \[ \left| h_{\varepsilon} \left( z_c + w_0 \varepsilon \ln
     \frac{1}{\varepsilon} \right) - h (z_c^+) \right| + \left|
     h_{\varepsilon}' \left( z_c + w_0 \varepsilon \ln \frac{1}{\varepsilon}
     \right) - h' (z_c^+) \right| \leqslant K \varepsilon \ln
     \frac{1}{\varepsilon} \]
\end{lemma}

This lemma implies that, when we are at distance $w_0 \varepsilon \ln
\frac{1}{\varepsilon}$ to the right of $z_c$ for some constant $w_0 > 0$
independent of $\varepsilon$, the function $h_{\varepsilon}$ and $h$ and their
derivatives are close. In particular, $h'_{\varepsilon} (z_c) = \frac{- a^2}{2
\varepsilon}$ is large when $\varepsilon$ is small, but $h_{\varepsilon}'
\left( z_c + w_0 \varepsilon \ln \frac{1}{\varepsilon} \right)$ is of size 1.
In other words, at $z_c + w_0 \varepsilon \ln \frac{1}{\varepsilon}$ the jump
has ended, and $h_{\varepsilon}, h'_{\varepsilon}$ will be for now on bounded
uniformly in $\varepsilon$. The choice of $w_0 \varepsilon \ln
\frac{1}{\varepsilon}$ is not necessarly optimal, it might be improved, but it
is enough here. We also compute the first order correction in $C^1$ between
$h_{\varepsilon}$ and $h$ to the right of $z_c$.

\begin{proof}
  We decompose, for $z > z_c, Z = \frac{z - z_c}{\varepsilon} > 0$, the
  solution of (\ref{41pb}) as
  \[ h_{\varepsilon} (z) = h (z) + (F (Z) + a) + \varepsilon G (Z), \]
  and we recall that $\frac{\alpha}{1 + \alpha} h (z) + \left( \frac{z}{1 +
  \alpha} - h (z) \right) h' (z) = 0.$ The function $h$ is discontinuous at $z
  = z_c$, but we focus here only on $z \in \left] z_c, z_c + w_0 \varepsilon
  \ln \frac{1}{\varepsilon} \right]$. We choose $F$ the solution of
  \[ \left\{\begin{array}{l}
       F'' (Z) - F (Z) F' (Z) = 0\\
       F (0) = 0, F' (0) = \frac{- a^2}{2},
     \end{array}\right. \]
  that is
  \[ F (Z) = \frac{a (e^{- a Z} - 1)}{1 + e^{- a Z}} . \]
  Remark that $F (+ \infty) = - a$, $a = \frac{h (z_c^-) - h (z_c^+)}{2}$ and
  $\frac{z_c}{1 + \alpha} = \frac{h (z_c^-) + h (z_c^+)}{2}$ by Lemma
  \ref{L42}. We also check that $G (0) = 0, G' (0) = - h' (z_c^+)$. Indeed,
  \[ h_{\varepsilon} (z^+_c) = h (z^+_c) + (F (0) + a) = \frac{z_c}{1 +
     \alpha} \]
  and
  \[ h_{\varepsilon}' (z^+_c) = h' (z_c^+) + \frac{1}{\varepsilon} F' (0) + G'
     (0) = \frac{- a^2}{2 \varepsilon} . \]
  Let us compute the equation satisfied by $G$. We have
  \[ \frac{\alpha}{1 + \alpha} h_{\varepsilon} = \frac{\alpha}{1 + \alpha} (h
     (z) + (F (Z) + a) + \varepsilon G (Z)), \]
  \[ \left( \frac{z}{1 + \alpha} - h_{\varepsilon} \right) = \left( \frac{z}{1
     + \alpha} - h (z) \right) - (F (Z) + a) - \varepsilon G (Z), \]
  \[ \partial_z h_{\varepsilon} = h' (z) + \frac{1}{\varepsilon} F' (Z) + G'
     (Z) \]
  so
  \begin{eqnarray*}
    &  & \left( \frac{z}{1 + \alpha} - h_{\varepsilon} \right) \partial_z
    h_{\varepsilon}\\
    & = & - \frac{1}{\varepsilon} F' (Z) F (Z) + \left( \frac{z}{1 + \alpha}
    - h (z) \right) h' (z)\\
    & + & \frac{1}{\varepsilon} F' (Z) \left( \frac{z}{1 + \alpha} - h (z) -
    a \right) + \left( \frac{z}{1 + \alpha} - h (z) \right) G' (Z)\\
    & - & (F (Z) + a) (h' (z) + G' (Z)) - G (Z) F' (Z)\\
    & - & \varepsilon G (Z) (h' (z) + G' (Z)) .
  \end{eqnarray*}
  Finally,
  \[ \varepsilon \partial_z^2 h_{\varepsilon} = \varepsilon h'' (z) +
     \frac{1}{\varepsilon} F'' (Z) + G'' (Z) . \]
  Using $F'' (Z) - F (Z) F' (Z) = 0$ and $\frac{\alpha}{1 + \alpha} h (z) +
  \left( \frac{z}{1 + \alpha} - h (z) \right) h' (z) = 0$ we infer that $G$
  satisfies on $Z > 0$ the equation
  \begin{eqnarray*}
    &  & G'' (Z) + \left( \frac{z}{1 + \alpha} - h (z) - (F (Z) + a) \right)
    G' (Z)\\
    & + & \left( - F' (Z) + \varepsilon \left( \frac{\alpha}{1 + \alpha} - h'
    (z) \right) \right) G (Z) - \varepsilon G (Z) G' (Z)\\
    & + & \frac{1}{\varepsilon} F' (Z) \left( \frac{z}{1 + \alpha} - h (z) -
    a \right) + \left( \frac{\alpha}{1 + \alpha} - h' (z) \right) (F (Z) +
    a)\\
    & + & \varepsilon h'' (z)\\
    & = & 0.
  \end{eqnarray*}
  We define the source part as
  \begin{eqnarray*}
    S & \assign & \frac{1}{\varepsilon} F' (Z) \left( \frac{z}{1 + \alpha} - h
    (z) - a \right) + \left( \frac{\alpha}{1 + \alpha} - h' (z) \right) (F (Z)
    + a)\\
    & + & \varepsilon h'' (z)
  \end{eqnarray*}
  and the operator on $G$ as
  \begin{eqnarray*}
    \mathcal{O} (G) & \assign & G'' (Z) + \left( \frac{z}{1 + \alpha} - h (z)
    - (F (Z) + a) \right) G' (Z)\\
    & + & \left( - F' (Z) + \varepsilon \left( \frac{\alpha}{1 + \alpha} - h'
    (z) \right) \right) G (Z) - \varepsilon G (Z) G' (Z),
  \end{eqnarray*}
  leading to the equation
  \[ \mathcal{O} (G) + S = 0. \]
  Let us estimate $S (Z)$ for $Z > 0$. We have
  \[ \frac{z_c}{1 + \alpha} - h (z^+_c) - a = 0, \]
  therefore
  \begin{eqnarray*}
    S (Z) & = & F' (Z) \left( \frac{Z}{1 + \alpha} - \frac{h (z_c +
    \varepsilon Z) - h (z_c)}{\varepsilon} \right)\\
    & + & \left( \frac{\alpha}{1 + \alpha} - h' (z_c + \varepsilon Z) \right)
    (F (Z) + a) + \varepsilon h'' (z_c + \varepsilon Z) .
  \end{eqnarray*}
  Take for now any $w_0 > 0$, independent of $\varepsilon$. Then, since
  \[ | F (Z) + a | + | F' (Z) | + | F'' (Z) | \leqslant K e^{- a Z}, \]
  we deduce that
  \[ | S (Z) | \leqslant K \left( \varepsilon + e^{- \frac{a}{2} Z} \right) \]
  for $Z \in \left] 0, w_0 \varepsilon \ln \frac{1}{\varepsilon} \right]$ for
  a constant $K > 0$ depending only on $w_0, \alpha, a$.
  
  \
  
  Let us now look at the coefficient in the operator $\mathcal{O} (G)$. We
  write it
  \[ \mathcal{O} (G) = G'' (Z) + A_1 (Z) G' (Z) + A_2 (Z) G (Z) - \varepsilon
     G (Z) G' (Z) \]
  with
  \[ A_1 (Z) \assign \frac{z_c + \varepsilon Z}{1 + \alpha} - h (z_c +
     \varepsilon Z) - (F (Z) + a) \]
  and
  \[ A_2 (Z) \assign - F' (Z) + \varepsilon \left( \frac{\alpha}{1 + \alpha} -
     h' (z_c + \varepsilon Z) \right) . \]
  In particular, $A_1$ and $A_2$ are bounded by constants independent of
  $\varepsilon$ if $\varepsilon < 1$. By the estimates on $S$, for any $Z_0 >
  0$, if $\varepsilon > 0$ is small enough depending on $Z_0$ (so that the
  nonlinear term $\varepsilon G (Z) G' (Z)$ can be neglected), there exists $K
  (Z_0) > 0$ such that
  \begin{equation}
    | G (Z) | + | G' (Z) | \leqslant K (Z_0) \label{420}
  \end{equation}
  for $Z \in [0, Z_0]$. This is because the equation satisfied by $G$ is,
  except for the term $- \varepsilon G G'$, linear with a bounded source term.
  Without this nonlinear term the solution would then be global, and taking
  $\varepsilon > 0$ small enough depending on $Z_0$, since $G (0), G' (0),
  A_1$ and $A_2$ are bounded uniformly in $\varepsilon$, the solution exists
  at least on $[0, Z_0]$ with a uniform estimates depending on $Z_0$.
  
  \
  
  Now, remark that $A_1 (Z) \rightarrow \frac{z_c}{1 + \alpha} - h (z^+_c)$
  and
  \[ \frac{A_2 (Z)}{\varepsilon} \rightarrow \left( \frac{\alpha}{1 + \alpha}
     - h' (z^+_c) \right) \]
  if $Z \geqslant w_0 \ln \frac{1}{\varepsilon}$ with $w_0$ large (such that
  $F' \left( w_0 \ln \frac{1}{\varepsilon} \right) \leqslant \varepsilon^2$
  for instance) when $\varepsilon \rightarrow 0$. We therefore write the
  equation on $G$ as
  \[ G'' (Z) + \left( \frac{z_c}{1 + \alpha} - h (z^+_c) \right) G' (Z) +
     \varepsilon \left( \frac{\alpha}{1 + \alpha} - h' (z^+_c) \right) G (Z) =
     S +\mathcal{R} (G) \]
  with
  \begin{eqnarray*}
    \mathcal{R} (G) & \assign & G' (Z) \left( \frac{z_c}{1 + \alpha} - h
    (z^+_c) - \left( \frac{z}{1 + \alpha} - h (z^+) - (F (Z) + a) \right)
    \right)\\
    & + & G (Z) \left( \varepsilon \left( \frac{\alpha}{1 + \alpha} - h'
    (z^+_c) \right) - \left( - F' (Z) + \varepsilon \left( \frac{\alpha}{1 +
    \alpha} - h' (z) \right) \right) \right)\\
    & + & \varepsilon G (Z) G' (Z) .
  \end{eqnarray*}
  We simplify
  \begin{eqnarray*}
    \mathcal{R} (G) & = & G' (Z) \left( \frac{- \varepsilon Z}{1 + \alpha} + h
    (z_c + \varepsilon Z) - h (z^+_c) + F (Z) + a \right)\\
    & + & G (Z) (F' (Z) + \varepsilon (h' (z) - h' (z^+_c)))\\
    & + & \varepsilon G (Z) G' (Z) .
  \end{eqnarray*}
  We define, to simplify the notations, $\lambda \assign \frac{z_c}{1 +
  \alpha} - h (z^+_c) > 0, \mu \assign \frac{\alpha}{1 + \alpha} - h' (z^+_c)
  > 0$ (since $h' (z^+_c) < 0$) so that the equation on $G$ can be written as
  \[ G'' (Z) + \lambda G' (Z) + \varepsilon \mu G (Z) = S +\mathcal{R} (G) .
  \]
  For $\varepsilon > 0$ small enough we have $\lambda^2 - 4 \varepsilon \mu >
  0$, and then we can write with
  \[ \lambda_{\pm} \assign \frac{- \lambda \pm \sqrt{\lambda^2 - 4 \varepsilon
     \mu}}{2} < 0, \]
  satisfying $\lambda_+ + \lambda_- = - \lambda$, $\lambda_+ \sim -
  \frac{\varepsilon \mu}{\lambda}, \lambda_- \sim - \lambda$ when $\varepsilon
  \rightarrow 0$, that (we recall that $G (0) = 0, G' (0) = - h' (z_c^+)$)
  \begin{eqnarray}
    G (Z) & = & \frac{- h' (z_c^+)}{\sqrt{\lambda^2 - 4 \varepsilon \mu}}
    (e^{\lambda_+ Z} - e^{\lambda_- Z}) \nonumber\\
    & + & e^{\lambda_+ Z} \int_0^Z e^{\lambda_- u} \left( \int_0^u e^{\lambda
    v} (S (v) +\mathcal{R} (G) (v)) d v \right) d u.  \label{4300}
  \end{eqnarray}
  Let us show that for $C_0 > 0$ large enough (independently of $\varepsilon$)
  and $\varepsilon$ small enough, we have
  \begin{equation}
    | G (Z) | + | G' (Z) | \leqslant C_0 \left( \varepsilon \ln
    \frac{1}{\varepsilon} + e^{- \frac{a}{2} Z} \right) \label{44}
  \end{equation}
  for $Z \in \left[ 0, w_0 \ln \frac{1}{\varepsilon} \right]$. This is true on
  $[0, Z_0]$ for some $Z_0 > 0$ by (\ref{420}). Now, if the result is not
  true, we denote $w_0 \ln \frac{1}{\varepsilon} \geqslant Z_c \geqslant Z_0$
  the first value such that this estimates becomes an equality. Then, on $[0,
  Z_c]$ we have
  \[ | S +\mathcal{R} (G) | \leqslant K (1 + C_0) \left( \varepsilon \ln
     \frac{1}{\varepsilon} + e^{- \frac{a}{2} Z} \right) \]
  and plotting this estimates in (\ref{4300}) leads to
  \[ | G (Z_c) | \leqslant K (1 + C_0 e^{- C \varepsilon Z_c}) \left(
     \varepsilon \ln \frac{1}{\varepsilon} + e^{- \frac{a}{2} Z_c} \right) .
  \]
  for some constants $K, C > 0$ independent of $\varepsilon$ and $C_0$. We can
  easily show a similar estimate on $G' (Z_c)$, up to an increase on $K_0$. We
  deduce that if $Z_0, C_0$ are large enough and $\varepsilon$ small enough,
  then
  \[ \frac{C_0}{2} \left( \varepsilon \ln \frac{1}{\varepsilon} + e^{-
     \frac{a}{2} Z_c} \right) > | G (Z_c) | + | G' (Z_c) | = C_0 \left(
     \varepsilon \ln \frac{1}{\varepsilon} + e^{- \frac{a}{2} Z_c} \right), \]
  which is a contradiction.
  
  \
  
  This completes the proof of (\ref{44}). Going back to $h_{\varepsilon} (z)
  = h (z) + (F (Z) + a) + \varepsilon G (Z),$we deduce that
  \[ \| h_{\varepsilon} (z) - h (z) - (F (Z) + a) \|_{L^{\infty} \left( \left]
     z_c, z_c + w_0 \varepsilon \ln \frac{1}{\varepsilon} \right] \right)}
     \rightarrow 0 \]
  when $\varepsilon \rightarrow 0$, and taking $w_0$ large enough, by Lemma
  \ref{SvD}

  \[ \left| h_{\varepsilon} \left( z_c + w_0 \varepsilon \ln
     \frac{1}{\varepsilon} \right) - h (z^+_c) \right| \leqslant \left| h
     \left( z_c + w_0 \varepsilon \ln \frac{1}{\varepsilon} \right) - h
     (z^+_c) \right| + K \varepsilon^2 \ln \frac{1}{\varepsilon} \leqslant K
     \varepsilon \ln \frac{1}{\varepsilon} . \]
  Finally,
  \[ h_{\varepsilon}' (z) = h' (z) + \frac{1}{\varepsilon} F' (Z) + G' (Z), \]
  leading to $\left| h_{\varepsilon}' (z) - \frac{1}{\varepsilon} F' (Z)
  \right| \leqslant C$ a constant independent of $\varepsilon$, and since if
  $w_0$ is large enough $\left| G' \left( w_0 \ln \frac{1}{\varepsilon}
  \right) \right| \leqslant K \varepsilon^{1 / 2}$ and $\frac{1}{\varepsilon}
  \left| F' \left( w_0 \ln \frac{1}{\varepsilon} \right) \right| \leqslant
  \frac{K}{\varepsilon} \varepsilon^{a w_0}$, we conclude the proof of this
  lemma by
  \[ \left| h'_{\varepsilon} \left( z_c + w_0 \varepsilon \ln
     \frac{1}{\varepsilon} \right) - h' (z^+_c) \right| \leqslant K
     \varepsilon \ln \frac{1}{\varepsilon} . \]
\end{proof}

By standard Cauchy theory, at fixed $z_c$ and $a$, $\varepsilon \rightarrow
h_{\varepsilon}$ is a smooth function. We conclude this subsection by some
estimates on $\partial_{\varepsilon} h_{\varepsilon}$.

\begin{lemma}
  \label{l445}For $\alpha \in] 0, 1 [, z_c > 0, \frac{z_c}{1 + \alpha} > a >
  \frac{\alpha}{1 + \alpha} z_c$, there exists $\varepsilon_0, C > 0$
  depending only on $\alpha, z_c, w_0$ such that, if $\varepsilon_0 >
  \varepsilon > 0$ and $h_{\varepsilon}$ is the solution of (\ref{41pb}) for
  these parameters, then
  \[ \varepsilon \rightarrow \left( h_{\varepsilon} \left( z_c + w_0
     \varepsilon \ln \frac{1}{\varepsilon} \right), \partial_z h_{\varepsilon}
     \left( z_c + w_0 \varepsilon \ln \frac{1}{\varepsilon} \right) \right)
     \in C^1 (] 0, \varepsilon_0 [, \mathbb{R}^2) \]
  with
  \[ \left| \partial_{\varepsilon} \left( h_{\varepsilon} \left( z_c + w_0
     \varepsilon \ln \frac{1}{\varepsilon} \right) \right) \right| + \left|
     \partial_{\varepsilon} \left( \partial_z h_{\varepsilon} \left( z_c + w_0
     \varepsilon \ln \frac{1}{\varepsilon} \right) \right) \right| \leqslant C
     \left( \ln \frac{1}{\varepsilon} \right)^2 . \]
  Furthermore, for $z \in \left] z_c, z_c + w_0 \varepsilon \ln
  \frac{1}{\varepsilon} \right]$ we have
  \[ | \varepsilon \partial_{\varepsilon} h_{\varepsilon} (z) | + \frac{|
     \varepsilon \partial_z \partial_{\varepsilon} h_{\varepsilon} (z) |}{\ln
     \frac{1}{\varepsilon}} \leqslant C e^{- \frac{a}{2} \left| \frac{z -
     z_c}{\varepsilon} \right|} . \]
\end{lemma}

\begin{proof}
  We recall that with $Z = \frac{z - z_c}{\varepsilon}$, we have
  \[ h_{\varepsilon} (z) = h (z) + F (Z) + a + \varepsilon G (Z, \varepsilon)
     . \]
  In the previous lemma, we did not write the dependency of $G$ in
  $\varepsilon$ since we did not differentiate with respect to it, but we do
  so here. We deduce that
  \[ \varepsilon \partial_{\varepsilon} h_{\varepsilon} = - Z F' (Z) +
     \varepsilon (G - Z \partial_Z G) + \varepsilon^2 \partial_{\varepsilon}
     G. \]
  With the explicit formula of $F$ and (\ref{44}) we check that for $z \in
  \left] z_c, z_c + w_0 \varepsilon \ln \frac{1}{\varepsilon} \right]$ we have
  \[ | - Z F' (Z) + \varepsilon (G - Z \partial_Z G) | \leqslant K e^{-
     \frac{a}{2} | Z |} \]
  where $K > 0$ is depending only on $\alpha, z_c, w_0$. Furthermore, from the
  proof of Lemma \ref{cigaretteandkerosine} we know that $G$ satisfies the
  equation
  \[ \partial_Z^2 G + \lambda \partial_Z G + \varepsilon \mu G = S
     +\mathcal{R} (G) \]
  with $\lambda = \frac{z_c}{1 + \alpha} - h (z^+_c) > 0, \mu =
  \frac{\alpha}{1 + \alpha} - h' (z^+_c) > 0$, hence
  \[ \partial_Z^2 \partial_{\varepsilon} G + \lambda \partial_Z
     \partial_{\varepsilon} G + \varepsilon \mu \partial_{\varepsilon} G =
     \partial_{\varepsilon} S + \partial_{\varepsilon} (\mathcal{R} (G)) - \mu
     G. \]
  We check that
  \[ | \varepsilon \partial_{\varepsilon} S (Z) - \mu G | \leqslant K \left(
     \varepsilon + e^{- \frac{a}{2} Z} \right) \]
  and by similar arguments as the proof of Lemma \ref{cigaretteandkerosine},
  we conclude that
  \[ | \varepsilon \partial_{\varepsilon} G (Z) | \leqslant K \left(
     \varepsilon + e^{- \frac{a}{2} Z} \right) \]
  for some constant $C > 0$ depending only on $\alpha, z_c, w_0$. Finally, we
  have
  \[ \partial_{\varepsilon} \left( h_{\varepsilon} \left( z_c + w_0
     \varepsilon \ln \frac{1}{\varepsilon} \right) \right) =
     \partial_{\varepsilon} h_{\varepsilon} \left( z_c + w_0 \varepsilon \ln
     \frac{1}{\varepsilon} \right) + w_0 \partial_{\varepsilon} \left(
     \varepsilon \ln \frac{1}{\varepsilon} \right) \partial_z h_{\varepsilon}
     \left( z_c + w_0 \varepsilon \ln \frac{1}{\varepsilon} \right), \]
  leading to
  \[ \left| \partial_{\varepsilon} \left( h_{\varepsilon} \left( z_c + w_0
     \varepsilon \ln \frac{1}{\varepsilon} \right) \right) \right| \leqslant K
     \ln \frac{1}{\varepsilon} . \]
  Similarly,
  \begin{eqnarray*}
    &  & \partial_{\varepsilon} \left( \partial_z h_{\varepsilon} \left( z_c
    + w_0 \varepsilon \ln \frac{1}{\varepsilon} \right) \right)\\
    & = & \partial_{\varepsilon} \partial_z h_{\varepsilon} \left( z_c + w_0
    \varepsilon \ln \frac{1}{\varepsilon} \right) + w_0 \partial_{\varepsilon}
    \left( \varepsilon \ln \frac{1}{\varepsilon} \right) \partial_z^2
    h_{\varepsilon} \left( z_c + w_0 \varepsilon \ln \frac{1}{\varepsilon}
    \right),
  \end{eqnarray*}
  and since $\frac{\alpha}{1 + \alpha} h_{\varepsilon} + \left( \frac{z}{1 +
  \alpha} - h_{\varepsilon} \right) \partial_z h_{\varepsilon} + \varepsilon
  \partial_z^2 h_{\varepsilon} = 0$, we have
  \begin{eqnarray*}
    &  & \partial_z^2 h_{\varepsilon} \left( z_c + w_0 \varepsilon \ln
    \frac{1}{\varepsilon} \right)\\
    & = & \frac{- 1}{\varepsilon} \left( \frac{\alpha}{1 + \alpha}
    (h_{\varepsilon} - h) \left( z_c + w_0 \varepsilon \ln
    \frac{1}{\varepsilon} \right) \right)\\
    & - & \frac{1}{\varepsilon} \left( \frac{z_c + w_0 \varepsilon \ln
    \frac{1}{\varepsilon}}{1 + \alpha} - h \left( z_c + w_0 \varepsilon \ln
    \frac{1}{\varepsilon} \right) \right) (\partial_z h_{\varepsilon} -
    \partial_z h) \left( z_c + w_0 \varepsilon \ln \frac{1}{\varepsilon}
    \right)\\
    & - & \frac{1}{\varepsilon} \partial_z h_{\varepsilon} \left( z_c + w_0
    \varepsilon \ln \frac{1}{\varepsilon} \right) (h - h_{\varepsilon}) \left(
    z_c + w_0 \varepsilon \ln \frac{1}{\varepsilon} \right),
  \end{eqnarray*}
  hence, by Lemma \ref{cigaretteandkerosine},
  \[ \left| \partial_z^2 h_{\varepsilon} \left( z_c + w_0 \varepsilon \ln
     \frac{1}{\varepsilon} \right) \right| \leqslant K \ln
     \frac{1}{\varepsilon} . \]
  This concludes the proof of this lemma.
\end{proof}

\subsubsection{Estimates in $\left[ z_c - w_0 \varepsilon \ln
\frac{1}{\varepsilon}, z_c \right[$}

\begin{lemma}
  \label{rubjck}For $\alpha \in] 0, 1 [, z_c > 0, \frac{z_c}{1 + \alpha} > a >
  \frac{\alpha}{1 + \alpha} z_c$, there exists $w_0 > 0$ depending on $\alpha,
  z_c, a$ such that $h_{\varepsilon}$, the solution of (\ref{41pb}), is well
  defined on $\left[ z_c - w_0 \varepsilon \ln \frac{1}{\varepsilon}, z_c
  \right[$, satisfies
  \[ \left\| h_{\varepsilon} (z) - h (z) + \frac{2 a}{1 + e^{- a \left(
     \frac{z - z_c}{\varepsilon} \right)}} \right\|_{L^{\infty} \left( \left[
     z_c - w_0 \varepsilon \ln \frac{1}{\varepsilon}, z_c \right[ \right)}
     \rightarrow 0 \]
  and
  \[ \left\| h'_{\varepsilon} (z) + \frac{2 a^2 e^{- a \left( \frac{z -
     z_c}{\varepsilon} \right)}}{\varepsilon \left( 1 + e^{- a \left( \frac{z
     - z_c}{\varepsilon} \right)} \right)^2} \right\|_{L^{\infty} \left(
     \left[ z_c - w_0 \varepsilon \ln \frac{1}{\varepsilon}, z_c \right[
     \right)} \leqslant C \]
  for some constant $C > 0$ depending only on $\alpha, z_c, a$ when
  $\varepsilon \rightarrow 0$, and also
  \[ \left| h_{\varepsilon} \left( z_c - w_0 \varepsilon \ln
     \frac{1}{\varepsilon} \right) - h (z_c^-) \right| + \left|
     h_{\varepsilon}' \left( z_c - w_0 \varepsilon \ln \frac{1}{\varepsilon}
     \right) - h' (z_c^-) \right| \leqslant C \varepsilon \ln
     \frac{1}{\varepsilon} . \]
\end{lemma}

\begin{proof}
  For $z < z_c$, keeping the notation $Z = \frac{z - z_c}{\varepsilon} < 0$,
  we decompose $h_{\varepsilon}$, solution of (\ref{41pb}) as
  \[ h_{\varepsilon} (z) = h (z) + (F (Z) - a) + \varepsilon G (Z) \]
  for the same function $F$ as in the proof of Lemma
  \ref{cigaretteandkerosine}, but another function $G$. We recall that the
  function $h$ is not continuous at $z_c$, and since we consider here $z <
  z_c$, it will have a different limit for $z - z_c < 0$ close to $0$. We take
  $G (0) = 0$ and $G' (0) = - h' (z_c^-)$ so that we match the conditions at
  $z_c$ of $h_{\varepsilon}$: $h_{\varepsilon} (z_c) = h (z_c^-) + F (0) - a =
  \frac{z_c}{1 + \alpha}$ and $h'_{\varepsilon} (z_c) = h' (z_c^-) +
  \frac{1}{\varepsilon} F' (0) + G' (0) = \frac{- a^2}{2 \varepsilon}$. As in
  the proof of Lemma \ref{cigaretteandkerosine}, we check that $G$ satisfies
  the equation
  \[ \mathcal{O} (G) + S = 0 \]
  with
  \begin{eqnarray*}
    S (Z) & \assign & \frac{1}{\varepsilon} F' (Z) \left( \frac{z}{1 + \alpha}
    - h (z) + a \right) + \left( \frac{\alpha}{1 + \alpha} - h' (z) \right) (F
    (Z) - a)\\
    & + & \varepsilon h'' (z)
  \end{eqnarray*}
  and
  \[ \mathcal{O} (G) = G'' (Z) + A_1 (Z) G' (Z) + A_2 (Z) G (Z) - \varepsilon
     G (Z) G' (Z) \]
  with
  \[ A_1 (Z) \assign \frac{z_c + \varepsilon Z}{1 + \alpha} - h (z_c +
     \varepsilon Z) - (F (Z) - a) \]
  and
  \[ A_2 (Z) \assign - F' (Z) + \varepsilon \left( \frac{\alpha}{1 + \alpha} -
     h' (z_c + \varepsilon Z) \right) . \]
  We now define
  \[ \tilde{G} (Z) = G (- Z), \]
  that satisfies the equation
  \[ \tilde{G}'' (Z) - A_1 (- Z) \tilde{G}' (Z) + A_2 (- Z) \tilde{G} (Z) +
     \varepsilon \tilde{G} (Z) \tilde{G}' (Z) = S (- Z) . \]
  We therefore consider $Z > 0$ in the rest of the proof. Now, remark that
  \[ - A_1 (- Z) \rightarrow - \left( \frac{z_c}{1 + \alpha} - h (z_c^-)
     \right) > 0 \]
  and
  \[ \frac{A_2 (- Z)}{\varepsilon} \rightarrow \left( \frac{\alpha}{1 +
     \alpha} - h' (z_c^-) \right) > 0 \]
  if $Z \geqslant w_0 \ln \frac{1}{\varepsilon}$ for $w_0$ large when
  $\varepsilon \rightarrow 0$. We therefore define $\lambda \assign - \left(
  \frac{z_c}{1 + \alpha} - h (z_c^-) \right) > 0$ and $\mu = \frac{\alpha}{1 +
  \alpha} - h' (z_c^-) > 0$, and we can complete the proof in a similar
  fashion as for Lemma \ref{cigaretteandkerosine}.
\end{proof}

\begin{lemma}
  \label{l465}For $\alpha \in] 0, 1 [, z_c > 0, \frac{z_c}{1 + \alpha} > a >
  \frac{\alpha}{1 + \alpha} z_c$, there exists $\varepsilon_0, C > 0$
  depending only on $a, z_c, w_0$ such that, if $\varepsilon_0 > \varepsilon >
  0$ and $h_{\varepsilon}$ is the solution of (\ref{41pb}) for these
  parameters, then
  \[ \varepsilon \rightarrow \left( h_{\varepsilon} \left( z_c - w_0
     \varepsilon \ln \frac{1}{\varepsilon} \right), h'_{\varepsilon} \left(
     z_c - w_0 \varepsilon \ln \frac{1}{\varepsilon} \right) \right) \in C^1
     (] 0, \varepsilon_0 [, \mathbb{R}^2) \]
  with
  \[ \left| \partial_{\varepsilon} \left( h_{\varepsilon} \left( z_c - w_0
     \varepsilon \ln \frac{1}{\varepsilon} \right) \right) \right| + \left|
     \partial_{\varepsilon} \left( h'_{\varepsilon} \left( z_c - w_0
     \varepsilon \ln \frac{1}{\varepsilon} \right) \right) \right| \leqslant C
     \left( \ln \frac{1}{\varepsilon} \right)^2 \]
  when $\varepsilon \rightarrow 0$. Furthermore, for $z \in \left[ z_c - w_0
  \varepsilon \ln \frac{1}{\varepsilon}, z_c \right[$ we have
  \[ | \varepsilon \partial_{\varepsilon} h_{\varepsilon} (z) | + \frac{|
     \varepsilon \partial_{\varepsilon} h_{\varepsilon} (z) |}{\ln
     \frac{1}{\varepsilon}} \leqslant C e^{- \frac{a}{2} \left| \frac{z -
     z_c}{\varepsilon} \right|} . \]
\end{lemma}

The proof of this result is similar to the proof of Lemma \ref{l445} and we
omit it.

\subsection{Profile far from $z_c$}\label{Timezone2}

\subsubsection{Profile on the right of $z_c$}

We start with an apriori estimate on solutions to the ODE problem.

\begin{lemma}
  \label{n6}For any $z_d > 0$, there exists $K > 0$ such that the solution to
  the problem
  \[ \left\{\begin{array}{l}
       \frac{\alpha}{1 + \alpha} h_{\varepsilon} + \left( \frac{z}{1 + \alpha}
       - h_{\varepsilon} \right) \partial_z h_{\varepsilon} + \varepsilon
       \partial_z^2 h_{\varepsilon} = 0\\
       \frac{z_d}{1 + \alpha} > h_{\varepsilon} (z_d) > 0, h_{\varepsilon}'
       (z_d) < 0
     \end{array}\right. \]
  for $\varepsilon > 0$ small enough (depending on $z_d, h_{\varepsilon}
  (z_d), h_{\varepsilon}' (z_d)$) is well defined on $[z_d, + \infty [$ and
  satisfies
  \[ | h_{\varepsilon} (z) | + z | h_{\varepsilon}' (z) | \leqslant
     \frac{K}{z^{\alpha}} \]
  for any $z \geqslant z_d$.
\end{lemma}

\begin{proof}
  Define
  \[ u = \frac{h_{\varepsilon}'}{h_{\varepsilon}} \]
  so that $\frac{h_{\varepsilon}''}{h_{\varepsilon}} = u' + u^2$ and
  $h_{\varepsilon} (z) = h_{\varepsilon} (z_d) e^{\int_{z_d}^z u (s) d s}$.
  Then, $u (z_d) = \frac{h_{\varepsilon}' (z_d)}{h_{\varepsilon} (z_d)} < 0$
  and
  \[ \varepsilon (u' + u^2) + \left( \frac{z}{1 + \alpha} - h_{\varepsilon}
     \right) u + \frac{\alpha}{1 + \alpha} = 0. \]
  We write it as
  \[ u' = - u^2 - \frac{1}{\varepsilon} \left( \left( \frac{z}{1 + \alpha} -
     h_{\varepsilon} \right) u + \frac{\alpha}{1 + \alpha} \right) . \]
  First, we have $u (z_d) < 0$, and we show that as long as $u$ exists, we
  have
  \[ u (z) < 0. \]
  Indeed, if $u (z) = 0$ for the first time at some point $z \geqslant z_d$,
  then $u' (z) = \frac{- \alpha}{(1 + \alpha) \varepsilon} < 0$, which is
  impossible. Using $h_{\varepsilon} (z) = h_{\varepsilon} (z_d)
  e^{\int_{z_d}^z u (s) d s}$, this implies that $h_{\varepsilon}$ is
  decreasing, and in particular $\frac{z}{1 + \alpha} - h_{\varepsilon} (z) >
  0$ for $z \geqslant z_d$.
  
  Also, for $\varepsilon > 0$ small enough, if $u (z_d) > \frac{- 1}{2
  \sqrt{\varepsilon}}$ say, then we always have $u (z) > \frac{- 1}{2
  \sqrt{\varepsilon}}$. This is because if at some point $u (z) = \frac{- 1}{2
  \sqrt{\varepsilon}}$, then $u' (z) > 0$, which is impossible. We deduce that
  $u$ is bounded, and therefore global.
  
  Using the same idea, we can show that $u (z) \leqslant \frac{- \lambda}{z}$
  for some small (but independent of $\varepsilon$ if $\varepsilon$ is small
  enough) constant $\lambda > 0$. In particular, $h_{\varepsilon} (z)
  z^{\lambda / 2} \rightarrow 0$ when $z \rightarrow + \infty$. Similarly, we
  can check that $u (z) \geqslant \frac{- 1}{\lambda z}$ if $\lambda$ is small
  enough.
  
  \
  
  Now define
  \[ v (z) = u (z) + \frac{\alpha}{z}, \]
  then
  \[ v' (z) + \frac{z}{(1 + \alpha) \varepsilon} v (z) = \frac{\alpha}{z^2} -
     u^2 + \frac{1}{\varepsilon} h_{\varepsilon} (z) u (z) . \]
  We have that
  \[ \left| \frac{\alpha}{z^2} - u^2 + \frac{1}{\varepsilon} h_{\varepsilon}
     (z) u (z) \right| \leqslant \frac{K}{\varepsilon z^{1 +
     \frac{\lambda}{2}}} \]
  and therefore, by a comparaison principle, on $z \geqslant z_d$ we have
  \[ | v (z) | \leqslant \frac{K}{z^{2 + \frac{\lambda}{2}}} . \]
  Using these estimates in the equation $h_{\varepsilon}' = u h_{\varepsilon}$
  completes the proof of the lemma.
\end{proof}

We recall that $h$ is solution of
\begin{equation}
  \left\{\begin{array}{l}
    \frac{\alpha}{1 + \alpha} h + \left( \frac{z}{1 + \alpha} - h \right)
    \partial_z h = 0\\
    h_0 (z) = \kappa_{\pm} (z_c, a) | z |^{- \alpha} (1 + o_{z \rightarrow \pm
    \infty} (1))
  \end{array}\right. \label{colorado}
\end{equation}
and is discontinuous at $z_c$.

\begin{lemma}
  \label{l488}The function $h_{\varepsilon}$, the solution of (\ref{41pb}),
  satisfies
  \[ \| (1 + | z |^{\alpha}) (h_{\varepsilon} (z) - h (z)) \|_{L^{\infty}
     \left( \left[ z_c + w_0 \varepsilon \ln \frac{1}{\varepsilon}, + \infty
     \right[ \right)} \rightarrow 0 \]
  when $\varepsilon \rightarrow 0$. Furhermore,
  \[ \lim_{z \rightarrow + \infty} z^{\alpha} h_{\varepsilon} (z) =:
     \kappa_{\varepsilon, +} (z_c, a) > 0 \]
  is well defined, and
  \[ | \kappa_{\varepsilon, +} (z_c, a) - \kappa_+ (z_c, a) | \rightarrow 0 \]
  when $\varepsilon \rightarrow 0$. Finally, for fixed values of $z_c$ and
  $a$, the function $\varepsilon \rightarrow \kappa_{\varepsilon, +} (z_c, a)$
  is smooth, and
  \[ | \partial_{\varepsilon} (\kappa_{\varepsilon, +} (z_c, a)) | \leqslant K
     \left( \ln \frac{1}{\varepsilon} \right)^2 \]
  for some constant $K$ independent of $\varepsilon$ if $\varepsilon$ is small
  enough. Furthermore,
  \[ \left| \partial_{\varepsilon} h_{\varepsilon} (z) -
     \frac{\partial_{\varepsilon} (\kappa_{\varepsilon, +} (z_c, a))}{| z
     |^{\alpha}} \right| \leqslant \frac{K \left( \ln \frac{1}{\varepsilon}
     \right)^2}{(1 + | z |)^{1 + 2 \alpha}} \]
  for $| z |$ large enough (uniformly in $\varepsilon$).
\end{lemma}

Remark that this does not imply that $\lim z^{\alpha} h_{\varepsilon} (z) =
\kappa_+ (z_c, a)$ when $z \rightarrow + \infty$, simply that their
differences goes to $0$ when $\varepsilon \rightarrow 0$.

\

\begin{proof}
  We introduce first a generic problem that we will use both to estimate
  $h_{\varepsilon}$ and $\partial_{\varepsilon} h_{\varepsilon}$.
  
  \
  
  We consider for now the problem
  \begin{equation}
    J_1 v + J_2 \partial_z v + \varepsilon \partial_z^2 v = S \label{NTO}
  \end{equation}
  for given functions $J_1, J_2, S$ of $z$, and initial values of $v$ at some
  point $z_d$, and with $J_2$ that does not cancel. We introduce the function
  $A$ defined by $A (z_d) = 1$ and
  \[ J_1 A + J_2 A' = 0, \]
  that is
  \[ A (z) = \exp \left( - \int_{z_d}^z \frac{J_1}{J_2} \right) > 0 \]
  Then, writing $v = A u$, we have
  \[ \varepsilon \partial_z^2 u + J_2 \partial_z u + \varepsilon \left( 2
     \frac{A'}{A} \partial_z u + \frac{A''}{A} u \right) = \frac{S}{A} . \]
  We continue, we introduce $B$ with $B (z_d) = 1$ and
  \[ \varepsilon B' = J_2 B, \]
  that is
  \[ B (z) = \exp \left( \frac{1}{\varepsilon} \int_{z_d}^z J_2 \right) > 0.
  \]
  We introduce for $\gamma > 0$ the quantity
  \[ \mathbb{B}_{\gamma} (z) \assign \frac{1}{B (z)} \int_{z_d}^z \frac{B
     (s)}{s^{\gamma}} d s, \]
  solution to the equation $\mathbb{B}'_{\gamma} (z) + \frac{J_2}{\varepsilon}
  \mathbb{B}_{\gamma} (z) = \frac{1}{z^{\gamma}}$ with $\mathbb{B}_{\gamma}
  (z_d) = 0$. If there exists a constant $C_0 > 0$ independent of
  $\varepsilon$ such that $\frac{1}{C_0} \geqslant \frac{J_2 (z)}{z} \geqslant
  C_0$, then by comparaison there exists $K > 0$ (depending only on $C_0$ and
  $\gamma$) such that
  \begin{equation}
    \mathbb{B}_{\gamma} (z) \leqslant \frac{K \varepsilon}{z^{\gamma + 1}} .
    \label{illuminate}
  \end{equation}
  We continue, we have
  \[ \varepsilon \partial_z^2 u + J_2 \partial_z u = \frac{\varepsilon
     \partial_z (\partial_z u B)}{B} \]
  and therefore
  \[ \partial_z (\partial_z u B) = B \left( \frac{S}{\varepsilon A} - 2
     \partial_z u \frac{A'}{A} + \frac{A''}{A} u \right) . \]
  Integrating between $z_d$ and $z$ leads to
  \begin{equation}
    \partial_z u = \frac{\partial_z u (z_d)}{B (z)} + \frac{1}{B (z)}
    \int_{z_d}^z B \left( \frac{S}{\varepsilon A} - 2 \partial_z u
    \frac{A'}{A} + \frac{A''}{A} u \right) . \label{paul}
  \end{equation}

  {\tmem{Step 1. Existence and properties of $\kappa_{\varepsilon, +} (z_c,
  a)$.}}
  
  \
  
  We take $z_d = z_c + w_0 \varepsilon \ln \frac{1}{\varepsilon}$, and we
  recall that $h_{\varepsilon}$ satisfies for $\varepsilon > 0$ small enough
  the equation
  \[ \frac{\alpha}{1 + \alpha} h_{\varepsilon} + \left( \frac{z}{1 + \alpha} -
     h_{\varepsilon} \right) \partial_z h_{\varepsilon} + \varepsilon
     \partial_z^2 h_{\varepsilon} = 0 \]
  and we have $\frac{z_d}{1 + \alpha} > h_{\varepsilon} (z_d) > 0,
  h_{\varepsilon}' (z_d) < 0$ since $h_{\varepsilon}$ solved (\ref{41pb}). We
  decompose $h_{\varepsilon} = h + g$ with $h$ solution of (\ref{colorado}),
  so that the equation satisfied by $g$ is
  \[ \varepsilon g'' + g' \left( \frac{z}{1 + \alpha} - h_{\varepsilon}
     \right) + g \left( \frac{\alpha}{1 + \alpha} - h' \right) + \varepsilon
     h'' = 0. \]
  This is equation (\ref{NTO}) with $J_1 = \frac{\alpha}{1 + \alpha} - h', J_2
  = \frac{z}{1 + \alpha} - h_{\varepsilon}$ and $\frac{S}{\varepsilon} = -
  h''$. Remark that for $z \geqslant z_d$, we have $J_2 (z) > 0$. Now, we have
  \[ \frac{J_1}{J_2} = \frac{\alpha - (1 + \alpha) h'}{z - (1 + \alpha)
     h_{\varepsilon}} = \frac{\alpha}{z} - \frac{(1 + \alpha) (z h' -
     h_{\varepsilon})}{z (z - (1 + \alpha) h_{\varepsilon})}, \]
  hence
  \[ A (z) = \exp \left( - \int_{z_d}^z \frac{J_1}{J_2} \right) = \left(
     \frac{z}{z_d} \right)^{- \alpha} \exp \left( - \int_{z_d}^z \frac{(1 +
     \alpha) (s h' - h_{\varepsilon})}{s (s - (1 + \alpha) h_{\varepsilon})} d
     s \right), \]
  and with
  \[ \left| \frac{(1 + \alpha) (s h' - h_{\varepsilon})}{s (s - (1 + \alpha)
     h_{\varepsilon})} \right| \leqslant \frac{K}{s^2} \]
  for $s \geqslant z_d$, we deduce that there exists $K > 0$ depending on
  $z_d, \alpha$ such that
  \[ \frac{1}{K} \leqslant z^{\alpha} (| A (z) | + | z A' (z) | + | z^2 A''
     (z) |) \leqslant K \]
  for $z \geqslant z_d$, and $z^{\alpha} A (z)$ converges when $z \rightarrow
  + \infty$ to a finite constant bounded uniformly in $\varepsilon$. With $g =
  A u$, we define $N (z) \assign \| u \|_{L^{\infty} ([z_d, z])} + \| z u'
  \|_{L^{\infty} ([z_d, z])}$. We have for $z_d \leqslant s \leqslant z$ that
  \[ \left| B \left( \frac{S}{\varepsilon A} - 2 \partial_z u \frac{A'}{A} +
     \frac{A''}{A} u \right) \right| (s) \leqslant \frac{B (s) (1 + N
     (z))}{s^2}, \]
  and by (\ref{illuminate}) we deduce that
  \[ \left| \frac{1}{B (z)} \int_{z_d}^z B \left( \frac{S}{\varepsilon A} - 2
     \partial_z u \frac{A'}{A} + \frac{A''}{A} u \right) \right| \leqslant
     \frac{K \varepsilon (1 + N (z))}{z^3} . \]
  Now, we have
  \[ \frac{1}{B (z)} = \exp \left( \frac{- 1}{\varepsilon} \int_{z_d}^z \left(
     \frac{s}{1 + \alpha} - h_{\varepsilon} (s) \right) d s \right) \leqslant
     \exp \left( \frac{- K (z - z_d)^2}{\varepsilon} \right) . \]
  Combining these estimates in (\ref{paul}) and the integral of (\ref{paul}),
  we deduce that
  \[ N (z) \leqslant C (\varepsilon + N (z_d)) \]
  for some constant $C > 0$ independent of $\varepsilon$ and for all $z
  \geqslant z_d$. Furthermore,
  \[ u (+ \infty) = u (z_d) + \int_{z_d}^{+ \infty} \frac{\partial_z u
     (z_d)}{B (s)} d s + \int_{z_d}^{+ \infty} \frac{1}{B (s)} \int_{z_d}^s B
     \left( \frac{S}{\varepsilon A} - 2 \partial_z u \frac{A'}{A} +
     \frac{A''}{A} u \right) \]
  is a finite quantity that satisfies $| u (+ \infty) | \leqslant K (| u (z_d)
  | + | \partial_z u (z_d) | + \varepsilon)$, and
  \begin{equation}
    | u (z) - u (+ \infty) | \leqslant \frac{K \varepsilon}{z^2}
    \label{mymind}
  \end{equation}
  for $z \geqslant z_d$. By Lemma \ref{cigaretteandkerosine}, we have $| g
  (z_d) | + | g' (z_d) | \leqslant K \varepsilon \ln \frac{1}{\varepsilon}$
  and thus $| u (z_d) | + | \partial_z u (z_d) | \leqslant K \varepsilon \ln
  \frac{1}{\varepsilon}$. We deduce that, defining
  \[ \kappa_{\varepsilon, +} (z_c, a) = \kappa_+ (z_c, a) + u (+ \infty)
     \lim_{z \rightarrow + \infty} (z^{\alpha} A (z)), \]
  we have
  \[ \lim_{z \rightarrow + \infty} z^{\alpha} h_{\varepsilon} (z) =
     \kappa_{\varepsilon, +} (z_c, a) . \]
  More precisely,
  \begin{eqnarray*}
    h_{\varepsilon} (z) & = & h (z) + A (z) u (z)\\
    & = & \frac{\kappa_{\varepsilon, +} (z_c, a)}{z^{\alpha}} + \left(
    \frac{z^{\alpha} A (z) - \lim_{x \rightarrow + \infty} (x^{\alpha} A
    (x))}{z^{\alpha}} \right) u (+ \infty)\\
    & + & \left( h (z) - \frac{\kappa_+ (z_c, a)}{z^{\alpha}} \right) + A (z)
    (u (z) - u (+ \infty)),
  \end{eqnarray*}
  and with the explicit definition of $A$, we have
  \[ \left| \frac{z^{\alpha} A (z) - \lim_{x \rightarrow + \infty} (x^{\alpha}
     A (x))}{z^{\alpha}} \right| \leqslant \frac{K}{z^{1 + 2 \alpha}}, \]
  and by Lemma \ref{Vegas} we have
  \[ \left| h (z) - \frac{\kappa_+ (z_c, a)}{z^{\alpha}} \right| \leqslant
     \frac{K}{z^{1 + 2 \alpha}} . \]
  With (\ref{mymind}), we conclude that for $z \geqslant z_d$,
  \[ | z^{\alpha} h_{\varepsilon} (z) - \kappa_{\varepsilon, +} (z_c, a) |
     \leqslant \frac{K}{z^{1 + \alpha}} \]
  and
  \[ | \kappa_{\varepsilon, +} (z_c, a) - \kappa_+ (z_c, a) | \rightarrow 0 \]
  when $\varepsilon \rightarrow 0$. Lemmas \ref{Vegas},
  \ref{cigaretteandkerosine} and $N (z) \leqslant C (\varepsilon + N (z_d))
  \leqslant K \varepsilon \ln \frac{1}{\varepsilon}$ also implies that for $z
  \geqslant z_d$,
  \[ | \partial_z h_{\varepsilon} (z) | \leqslant \frac{K}{z^{1 + \alpha}}, \]
  as well as
  \[ | \partial^2_z h_{\varepsilon} (z) | \leqslant \frac{K \ln
     \frac{1}{\varepsilon}}{z^{2 + \alpha}} . \]
  This last estimate can be improved (we can remove the $\ln
  \frac{1}{\varepsilon}$) but it is not needed here, we will only need
  $\varepsilon \partial_{\varepsilon} h_{\varepsilon}$ to be small and not
  $\partial_{\varepsilon} h_{\varepsilon}$ itself.
  
  \
  
  {\tmem{Step 2. Differentiation with respect to $\varepsilon$ at fixed
  $z_d$.}}
  
  \
  
  We consider here $\mathfrak{h}_{\varepsilon}$ the solution to the problem
  \[ \frac{\alpha}{1 + \alpha} \mathfrak{h}_{\varepsilon} + \left( \frac{z}{1
     + \alpha} -\mathfrak{h}_{\varepsilon} \right) \partial_z
     \mathfrak{h}_{\varepsilon} + \varepsilon \partial_z^2
     \mathfrak{h}_{\varepsilon} = 0 \]
  with $\mathfrak{h}_{\varepsilon} (z_d), \partial_z
  \mathfrak{h}_{\varepsilon} (z_d)$ given satisfying $\frac{z_d}{1 + \alpha}
  >\mathfrak{h}_{\varepsilon} (z_d) > 0, \mathfrak{h}_{\varepsilon}' (z_d) <
  0$, and they are, with $z_d$, independent of $\varepsilon$. We have as
  previously that for $z \geqslant z_d$,
  \begin{equation}
    | \mathfrak{h}_{\varepsilon} (z) | + z | \partial_z
    \mathfrak{h}_{\varepsilon} (z) | \leqslant \frac{K}{z^{\alpha}}, |
    \partial^2_z \mathfrak{h}_{\varepsilon} (z) | \leqslant \frac{K \ln
    \frac{1}{\varepsilon}}{z^{2 + \alpha}} . \label{surrender}
  \end{equation}
  We introduce this new notation since we want to differentiate
  $h_{\varepsilon}$ with respect to $\varepsilon$, but its dependency on
  $\varepsilon$ comes from the $\varepsilon$ in the equation but also from
  $z_c$ and the value of $h_{\varepsilon}$ here. For
  $\mathfrak{h}_{\varepsilon}$, the dependency on $\varepsilon$ comes only
  from the $\varepsilon$ in front of $\partial_z^2 \mathfrak{h}_{\varepsilon}$
  in the equation. By standard Cauchy-Theory, $\varepsilon \rightarrow
  \mathfrak{h}_{\varepsilon}$ is differentiable, and $v =
  \partial_{\varepsilon} \mathfrak{h}_{\varepsilon}$ satisfies the problem
  \[ \left\{\begin{array}{l}
       \left( \frac{\alpha}{1 + \alpha} - \partial_z
       \mathfrak{h}_{\varepsilon} \right) v + \left( \frac{z}{1 + \alpha}
       -\mathfrak{h}_{\varepsilon} \right) \partial_z v + \varepsilon
       \partial_z^2 v = - \partial_z^2 \mathfrak{h}_{\varepsilon}\\
       v (z_d) = v' (z_d) = 0.
     \end{array}\right. \]
  This is equation (\ref{NTO}) with $J_1 = \frac{\alpha}{1 + \alpha} -
  \partial_z \mathfrak{h}_{\varepsilon}, J_2 = \frac{z}{1 + \alpha}
  -\mathfrak{h}_{\varepsilon}$ and $\frac{S}{\varepsilon} = \frac{-
  1}{\varepsilon} \partial_z^2 \mathfrak{h}_{\varepsilon} .$ Following a
  similar proof as step 1, we check that, with $\mathfrak{h}_{\varepsilon} = A
  u$, we have that $u (z)$ converges to a finite limite $u (+ \infty)$, with
  $| u (+ \infty) | \leqslant K \ln \frac{1}{\varepsilon}$, and that
  \[ | u (z) - u (+ \infty) | \leqslant \frac{K \ln
     \frac{1}{\varepsilon}}{z^2} \]
  for $z \geqslant z_d$. We also check, as in step 1, that
  \[ | z^{\alpha} v (z) - k_0 | \leqslant \frac{K \ln
     \frac{1}{\varepsilon}}{z^{1 + \alpha}} \]
  and
  \[ | \partial_z v | \leqslant \frac{K \ln \frac{1}{\varepsilon}}{z^{1 +
     \alpha}} \]
  for some $k_0$ depending on $\varepsilon$ and $K > 0$.
  
  \
  
  {\tmem{Step 3. Differentiation with respect to $z_d$.}}
  
  \
  
  We consider here $\mathfrak{h}_{\varepsilon}$ the solution to the problem
  \[ \frac{\alpha}{1 + \alpha} \mathfrak{h}_{\varepsilon} + \left( \frac{z}{1
     + \alpha} -\mathfrak{h}_{\varepsilon} \right) \partial_z
     \mathfrak{h}_{\varepsilon} + \varepsilon \partial_z^2
     \mathfrak{h}_{\varepsilon} = 0 \]
  with $\mathfrak{h}_{\varepsilon} (z_d) = a, \partial_z
  \mathfrak{h}_{\varepsilon} (z_d) = b$, where $\frac{z_d}{1 + \alpha} > a >
  0, b < 0$ are independent of $z_d$ and verify $\left| \frac{\alpha}{1 +
  \alpha} a + \left( \frac{z_d}{1 + \alpha} - a \right) b \right| \leqslant K
  \varepsilon \ln \frac{1}{\varepsilon}$ for some $K > 0$ independent of
  $\varepsilon$. As previously, estimate (\ref{surrender}) holds. We want to
  compute $v = \partial_{z_d} \mathfrak{h}_{\varepsilon}$. It is solution of
  \[ \left\{\begin{array}{l}
       \left( \frac{\alpha}{1 + \alpha} - \partial_z
       \mathfrak{h}_{\varepsilon} \right) v + \left( \frac{z}{1 + \alpha}
       -\mathfrak{h}_{\varepsilon} \right) \partial_z v + \varepsilon
       \partial_z^2 v = 0\\
       v (z_d) = - \partial_z \mathfrak{h}_{\varepsilon} (z_d), v' (z_d) = -
       \partial_z^2 \mathfrak{h}_{\varepsilon} (z_d) .
     \end{array}\right. \]
  We have $- \partial_z \mathfrak{h}_{\varepsilon} (z_d) = - b$ and since
  $\left( \frac{\alpha}{1 + \alpha} \mathfrak{h}_{\varepsilon} + \left(
  \frac{z}{1 + \alpha} -\mathfrak{h}_{\varepsilon} \right) \partial_z
  \mathfrak{h}_{\varepsilon} + \varepsilon \partial_z^2
  \mathfrak{h}_{\varepsilon} \right) (z_d) = 0$, we have
  \[ - \partial_z^2 \mathfrak{h}_{\varepsilon} (z_d) = \frac{1}{\varepsilon}
     \left( \frac{\alpha}{1 + \alpha} a + \left( \frac{z_d}{1 + \alpha} - a
     \right) b \right) \]
  which is bounded by $K \ln \frac{1}{\varepsilon}$ with $K > 0$ independent
  of $\varepsilon$. As in the previous case, we check that
  \[ | z^{\alpha} v (z) -k_0 | \leqslant \frac{K \ln
     \frac{1}{\varepsilon}}{z^{1 + \alpha}} \]
  for some $k_0, K > 0$ and
  \[ | \partial_z v | \leqslant \frac{K \ln \frac{1}{\varepsilon}}{z^{1 +
     \alpha}} \]

  {\tmem{Step 4. Differentiation with respect to $\mathfrak{h}_{\varepsilon}
  (z_d)$ and $\partial_z \mathfrak{h}_{\varepsilon} (z_d)$.}}
  
  \
  
  We consider here $\mathfrak{h}_{\varepsilon}$ the solution to the problem
  \[ \frac{\alpha}{1 + \alpha} \mathfrak{h}_{\varepsilon} + \left( \frac{z}{1
     + \alpha} -\mathfrak{h}_{\varepsilon} \right) \partial_z
     \mathfrak{h}_{\varepsilon} + \varepsilon \partial_z^2
     \mathfrak{h}_{\varepsilon} = 0 \]
  with $\mathfrak{h}_{\varepsilon} (z_d) = a, \partial_z
  \mathfrak{h}_{\varepsilon} (z_d) = b$, where $\frac{z_d}{1 + \alpha} > a >
  0, b < 0$. Then, $v = \partial_a \mathfrak{h}_{\varepsilon}$ satisfies
  \[ \left\{\begin{array}{l}
       \left( \frac{\alpha}{1 + \alpha} - \partial_z
       \mathfrak{h}_{\varepsilon} \right) v + \left( \frac{z}{1 + \alpha}
       -\mathfrak{h}_{\varepsilon} \right) \partial_z v + \varepsilon
       \partial_z^2 v = 0\\
       v (z_d) = 1, v' (z_d) = 0.
     \end{array}\right. \]
  This is similar to the previous steps, and we also check that $\partial_b
  \mathfrak{h}_{\varepsilon}$ can be estimates similarly
  
  \
  
  {\tmem{Step 5. Conclusion.}}
  
  \
  
  The function $h_{\varepsilon}$ is solution to $\frac{\alpha}{1 + \alpha}
  h_{\varepsilon} + \left( \frac{z}{1 + \alpha} - h_{\varepsilon} \right)
  \partial_z h_{\varepsilon} + \varepsilon \partial_z^2 h_{\varepsilon} = 0$
  with the initial condition at $z_d = z_c + w_0 \varepsilon \ln
  \frac{1}{\varepsilon}$ satisfying (by Lemmas \ref{cigaretteandkerosine} and
  \ref{l445})
  \[ | h_{\varepsilon} (z_d) - h (z_c^+) | + | h_{\varepsilon}' (z_d) - h'
     (z_c^+) | \leqslant K \varepsilon \ln \frac{1}{\varepsilon} \]
  and
  \[ | \partial_{\varepsilon} h_{\varepsilon} (z_d) | + |
     \partial_{\varepsilon} h'_{\varepsilon} (z_d) | \leqslant K \left( \ln
     \frac{1}{\varepsilon} \right)^2 . \]
  Therefore, $\partial_{\varepsilon} h_{\varepsilon}$ can be written as a sum
  of the functions $v$ of steps 2 to 4, and since $| \partial_{\varepsilon}
  z_d | \leqslant w_0 \ln \frac{1}{\varepsilon}$, this concludes the proof of
  this lemma.
\end{proof}

\subsubsection{Profile on the left of $z_c$}

\begin{lemma}
  \label{rubick}The function $h_{\varepsilon}$, the solution of (\ref{41pb}),
  is well defined on $\left] - \infty, z_c - w_0 \varepsilon \ln
  \frac{1}{\varepsilon} \right]$ and satisfies
  \[ \| (1 + | z |^{\alpha}) (h_{\varepsilon} (z) - h (z)) \|_{L^{\infty}
     \left( \left] - \infty, z_c - w_0 \varepsilon \ln \frac{1}{\varepsilon}
     \right] \right)} \rightarrow 0 \]
  when $\varepsilon \rightarrow 0$. Furhermore,
  \[ \lim_{z \rightarrow - \infty} z^{\alpha} h_{\varepsilon} (z) =:
     \kappa_{\varepsilon, -} (z_c, a) > 0 \]
  is well defined, and
  \[ | \kappa_{\varepsilon, -} (z_c, a) - \kappa_- (z_c, a) | \rightarrow 0 \]
  when $\varepsilon \rightarrow 0$. Finally, for fixed values of $z_c$ and
  $a$, the function $\varepsilon \rightarrow \kappa_{\varepsilon, -} (z_c, a)$
  is smooth, and
  \[ | \partial_{\varepsilon} (\kappa_{\varepsilon, -} (z_c, a)) | \leqslant K
     \left( \ln \frac{1}{\varepsilon} \right)^2 \]
  For a constant $K > 0$ independent of $\varepsilon$ if $\varepsilon$ is
  small enough. Furthermore,
  \[ \left| \partial_{\varepsilon} h_{\varepsilon} (z) -
     \frac{\partial_{\varepsilon} (\kappa_{\varepsilon, -} (z_c, a))}{| z
     |^{\alpha}} \right| \leqslant \frac{K \left( \ln \frac{1}{\varepsilon}
     \right)^2}{(1 + | z |)^{1 + 2 \alpha}} \]
  for $| z |$ large enough.
\end{lemma}

The proof is similar to the one of Lemma \ref{l488} and we omit it.

\subsection{End of the proof of Proposition
\ref{profileODEd}}\label{Timezone3}

Take $\kappa_+, \kappa_- > 0, \alpha \in] 0, 1 [$. By Lemma \ref{L42}, we
choose $z_c, a > 0$ such that
\[ \kappa_+ (z_c, a) = \kappa_+, \kappa_- (z_c, a) = \kappa_- . \]
We infer that for $\varepsilon$ small enough, we can take $z_c (\varepsilon),
a (\varepsilon) > 0$ such that
\[ \kappa_{\varepsilon, +} (z_c (\varepsilon), a (\varepsilon)) = \kappa_+,
   \kappa_{\varepsilon, -} (z_c (\varepsilon), a (\varepsilon)) = \kappa_- \]
with
\[ | z_c (\varepsilon) - z_{\kappa} | + | a (\varepsilon) - a_{\kappa} |
   \rightarrow 0 \]
when $\varepsilon \rightarrow 0$, and that this choice is unique, and
determined $\mathbbmss{h}_{\varepsilon}$. This is a consequence of the implicit
function theorem on the function
\[ \mathfrak{K} (\varepsilon, z_c, a) \assign (\kappa_{+, \varepsilon} (z_c,
   a) - \kappa_+, \kappa_{-, \varepsilon} (z_c, a) - \kappa_-) . \]
Indeed, by Lemmas \ref{L42}, \ref{l488} and \ref{rubick}, we have
$\mathfrak{K} (0, z_c, a) = 0$ and the jacobian at $\varepsilon = 0$ is
invertible. By Lemma \ref{l488} and \ref{rubick} we have the estimates
\begin{equation}
  | \partial_{\varepsilon} z_c (\varepsilon) | + | \partial_{\varepsilon} a
  (\varepsilon) | \leqslant K \left( \ln \frac{1}{\varepsilon} \right)^2
  \label{Pulverturm}
\end{equation}
The other properties in Proposition \ref{profileODEd} are a consequence of
Lemmas \ref{cigaretteandkerosine}, \ref{rubjck}, \ref{l488} and \ref{rubick}.

\subsection{Properties of $\partial_{\varepsilon} \mathbbmss{h}_{\varepsilon}$}

We recall that the function $\mathbbmss{h}_{\varepsilon}$ is a solution
$h_{\varepsilon}$ of the previous subsection, with a particular choice of $z_c
(\varepsilon), a (\varepsilon)$ such that the limits at $\pm \infty$ of $| z
|^{\alpha} \mathbbmss{h}_{\varepsilon} (z)$ are $\kappa_{\pm}$ respectively,
quantities independent of $\varepsilon$. The function
$\mathbbmss{h}_{\varepsilon}$ depends on $\varepsilon$ by $h_{\varepsilon}$ as
above, but also through $z_c (\varepsilon)$ and $a (\varepsilon)$.

\begin{lemma}
  \label{reboot}The function $\partial_{\varepsilon} \mathbbmss{h}_{\varepsilon}$
  satisfies
  \[ | \partial_{\varepsilon} \mathbbmss{h}_{\varepsilon} (z) | \leqslant \frac{K
     \left( \ln \frac{1}{\varepsilon} \right)^2}{(1 + | z |)^{1 + 2 \alpha}}
  \]
  for $| z - z_c (\varepsilon) | \geqslant 1$ and a constant $K$ independent
  of $\varepsilon$. For $z \in [z_c (\varepsilon) - 1, z_c (\varepsilon) + 1]$
  we have
  \[ | \varepsilon \partial_{\varepsilon} \mathbbmss{h}_{\varepsilon} (z) |
     \leqslant K e^{- \frac{a}{2} \left| \frac{z - z_c
     (\varepsilon)}{\varepsilon} \right|} . \]
  Finally, $\partial_{\varepsilon} \mathbbmss{h}_{\varepsilon} \in L^1
  (\mathbb{R})$ and
  \[ \int_{\mathbb{R}} \partial_{\varepsilon} \mathbbmss{h}_{\varepsilon} = 0. \]
\end{lemma}

\begin{proof}
  By Lemma \ref{l488}, we check that for $z \geqslant z_c (\varepsilon) + 1$,
  we have
  \[ \left| \partial_{\varepsilon} \mathbbmss{h}_{\varepsilon} (z) -
     \frac{\partial_{\varepsilon} (\kappa_{+, \varepsilon} (z_c (\varepsilon),
     a (\varepsilon)))}{| z |^{\alpha}} \right| \leqslant \frac{K \left( \ln
     \frac{1}{\varepsilon} \right)^2}{(1 + | z |)^{1 + 2 \alpha}} \]
  and
  \[ | \partial_{\varepsilon} \partial_z \mathbbmss{h}_{\varepsilon} (z) |
     \leqslant \frac{K \left( \ln \frac{1}{\varepsilon} \right)^2}{(1 + | z
     |)^{1 + \alpha}} \]
  but $\kappa_{+, \varepsilon} (z_c (\varepsilon), a (\varepsilon)) =
  \kappa_+$ which is independent of $\varepsilon$, hence
  $\partial_{\varepsilon} (\kappa_{+, \varepsilon} (z_c (\varepsilon), a
  (\varepsilon))) = 0$. Now, on $[z_c (\varepsilon), z_c (\varepsilon) + 1]$,
  the estimate is a consequence of Lemma \ref{l445}. For $z \leqslant z_c
  (\varepsilon)$ the proof is similar. The decay at infinity of \
  $\partial_{\varepsilon} \mathbbmss{h}_{\varepsilon}$ implies in particular that
  it is in $L^1 (\mathbb{R})$.
  
  Now, $\mathbbmss{h}_{\varepsilon}$ satisfies $\frac{\alpha}{1 + \alpha}
  \mathbbmss{h}_{\varepsilon} + \frac{z}{1 + \alpha} \partial_z
  \mathbbmss{h}_{\varepsilon} -\mathbbmss{h}_{\varepsilon} \partial_z
  \mathbbmss{h}_{\varepsilon} + \varepsilon \partial_z^2 \mathbbmss{h}_{\varepsilon}
  = 0,$ and integrating between $- x$ and $x$ for some large $x > 0$ leads to
  \[ \frac{\alpha - 1}{1 + \alpha} \int_{- x}^x \mathbbmss{h}_{\varepsilon} +
     \left[ \frac{z}{1 + \alpha} \mathbbmss{h}_{\varepsilon} - \frac{1}{2}
     \mathbbmss{h}_{\varepsilon}^2 + \partial_z \mathbbmss{h}_{\varepsilon}
     \right]_{- x}^x = 0. \]
  Differentiating with respect to $\varepsilon$ leads to
  \[ \frac{\alpha - 1}{1 + \alpha} \int_{- x}^x \partial_{\varepsilon}
     \mathbbmss{h}_{\varepsilon} + \left[ \frac{z}{1 + \alpha}
     \partial_{\varepsilon} \mathbbmss{h}_{\varepsilon} - \partial_{\varepsilon}
     \mathbbmss{h}_{\varepsilon} \mathbbmss{h}_{\varepsilon} + \partial_z
     \partial_{\varepsilon} \mathbbmss{h}_{\varepsilon} \right]_{- x}^x = 0, \]
  and going to the limit $x \rightarrow + \infty$, we check with $|
  \partial_{\varepsilon} \mathbbmss{h}_{\varepsilon} (z) | \leqslant \frac{K
  \left( \ln \frac{1}{\varepsilon} \right)^2}{(1 + | z |)^{1 + 2 \alpha}}$
  that
  \[ \frac{\alpha - 1}{1 + \alpha} \int_{\mathbb{R}} \partial_{\varepsilon}
     \mathbbmss{h}_{\varepsilon} = 0. \]
\end{proof}

\section{Stability of $\mathbbmss{h}_{\varepsilon}$}\label{Amywiles}

This section is devoted to the proof of Theorem \ref{prop111}. For the viscous
Burgers equation $\partial_t u - \partial_x^2 u + u \partial_x u = 0$ and
$\varepsilon (t) = t^{\frac{\alpha - 1}{\alpha + 1}}$, we introduce now the
rescaling
\[ y = \frac{x - t^{\frac{1}{1 + \alpha}} z_c}{\varepsilon} = t^{\frac{2 -
   \alpha}{1 + \alpha}} (z - z_c) \]
and the rescaled function
\[ H (y, t) = t^{\frac{\alpha}{1 + \alpha}} u \left( (z_c + \varepsilon (t) y)
   t^{\frac{1}{1 + \alpha}}, t \right) . \]
This scale is the right one to understand the profile near the discontinuity
point $z_c$ ($y = 0$ correspond to $z = z_c$). The previous scaling (in $z$),
where the profile $h_{\varepsilon}$ was constructed, was
\[ h (z, \varepsilon (t)) = t^{\frac{\alpha}{1 + \alpha}} u \left( z
   t^{\frac{1}{1 + \alpha}}, t \right), \]
and they are connected by
\[ H (y, t) = h (z_c (t) + \varepsilon (t) y, \varepsilon (t)) . \]
In particular, we define
\[ H_{\varepsilon} (y, t) \assign \mathbbmss{h}_{\varepsilon} (z_c (t) +
   \varepsilon (t) y) . \]
We recall that $\mathbbmss{h}_{\varepsilon}$ satisfies
\[ \frac{\alpha}{1 + \alpha} \mathbbmss{h}_{\varepsilon} + \frac{z}{1 + \alpha}
   \partial_z \mathbbmss{h}_{\varepsilon} -\mathbbmss{h}_{\varepsilon} \partial_z
   \mathbbmss{h}_{\varepsilon} + \varepsilon \partial_z^2
   \mathbbmss{h}_{\varepsilon} = 0, \]
therefore $H_{\varepsilon}$ satisfies
\[ \frac{\alpha}{1 + \alpha} H_{\varepsilon} + \frac{t^{\frac{1 - \alpha}{1 +
   \alpha}} z_c (t) + y}{1 + \alpha} \partial_y H_{\varepsilon} - t^{\frac{1 -
   \alpha}{1 + \alpha}} H_{\varepsilon} \partial_y H_{\varepsilon} +
   t^{\frac{1 - \alpha}{1 + \alpha}} \partial_y^2 H_{\varepsilon} = 0, \]
that is
\[ t^{- \frac{2 \alpha}{1 + \alpha}} (- \partial_y^2 H_{\varepsilon} +
   H_{\varepsilon} \partial_y H_{\varepsilon}) - \frac{\alpha}{1 + \alpha}
   t^{- 1} H_{\varepsilon} - t^{- 1} \frac{\left( t^{\frac{1 - \alpha}{1 +
   \alpha}} z_c (t) + y \right)}{1 + \alpha} \partial_y H_{\varepsilon} = 0.
\]
Now, we have
\[ \partial_y H = t^{\frac{2 \alpha}{1 + \alpha}} \partial_x u, \partial_y^2 H
   = t^{\frac{3 \alpha}{1 + \alpha}} \partial^2_x u \]
and
\[ \partial_t H = t^{\frac{\alpha}{1 + \alpha}} \partial_t u + \frac{\alpha}{1
   + \alpha} t^{- 1} H + \frac{t^{- 1} \left( \alpha y + t^{\frac{1 -
   \alpha}{1 + \alpha}} z_c \right)}{1 + \alpha} \partial_y H + \partial_t z_c
   t^{\frac{1 - \alpha}{1 + \alpha}} \partial_y H, \]
therefore the equation on $H$ (that is the rescaled viscous Burgers equation
in this new scaling) is
\begin{eqnarray*}
  &  & \partial_t H - \frac{\alpha}{1 + \alpha} t^{- 1} \partial_y (y H) -
  \frac{t^{- 1 + \frac{1 - \alpha}{1 + \alpha}} z_c}{1 + \alpha} \partial_y
  H\\
  & - & \partial_t z_c t^{\frac{1 - \alpha}{1 + \alpha}} \partial_y H + t^{-
  \frac{2 \alpha}{1 + \alpha}} (- \partial_y^2 H + H \partial_y H)\\
  & = & 0.
\end{eqnarray*}
We now decompose $H = H_{\varepsilon} + f$. Then $f$ satisfies the equation
\begin{eqnarray*}
  &  & \partial_t f - \frac{\alpha}{1 + \alpha} t^{- 1} \partial_y (y f)\\
  & + & t^{- \frac{2 \alpha}{1 + \alpha}} \left( - \partial_y^2 f + f
  \partial_y H_{\varepsilon} + \left( H_{\varepsilon} - \frac{z_c}{1 + \alpha}
  \right) \partial_y f + f \partial_y f \right) - \partial_t z_c t^{\frac{1 -
  \alpha}{1 + \alpha}} \partial_y f\\
  & + & \partial_t H_{\varepsilon} - \partial_t z_c t^{\frac{1 - \alpha}{1 +
  \alpha}} \partial_y H_{\varepsilon} + \frac{1 - \alpha}{1 + \alpha} t^{- 1}
  y \partial_y H_{\varepsilon}\\
  & = & 0.
\end{eqnarray*}
We check that
\[ S \assign \partial_t H_{\varepsilon} - \partial_t z_c t^{\frac{1 -
   \alpha}{1 + \alpha}} \partial_y H_{\varepsilon} + \frac{1 - \alpha}{1 +
   \alpha} t^{- 1} y \partial_y H_{\varepsilon} = \partial_t \varepsilon
   \partial_{\varepsilon} \mathbbmss{h}_{\varepsilon} (z_c + \varepsilon y) . \]
Remark that the problem takes the form
\begin{equation}
  t^{\frac{2 \alpha}{1 + \alpha}} \partial_t f + \partial_y \left( -
  \partial_y f + \frac{- \alpha \varepsilon}{1 + \alpha} y f + \left(
  H_{\varepsilon} - \frac{z_c}{1 + \alpha} \right) f + \frac{f^2}{2} - t
  \partial_t z_c f \right) + S = 0. \label{will}
\end{equation}
From Lemma \ref{reboot}, $S \in L^1 (\mathbb{R})$ and $\int_{\mathbb{R}} S =
0$. We therefore write $S = \partial_y \tilde{S}$ with
\[ \tilde{S} \assign \int_{- \infty}^y S. \]
However, the perturbation $f$ can have a mass. To deal with it, we introduce
an additional term $\mathfrak{u}_M$ in the profile, and we will decompose $f
=\mathfrak{u}_M + \partial_x g$, where all the mass is in $\mathfrak{u}_M$.

\subsection{Definition and properties of $\mathfrak{u}_M$}\label{Voyager}

This subsection is devoted to the proof of the following result.

\begin{lemma}
  Given $M \in \mathbb{R}, a > 0$, the problem
  \[ \left\{\begin{array}{l}
       - \partial_x \mathfrak{u}+ F\mathfrak{u}+ \frac{\mathfrak{u}^2}{2} =
       0\\
       \int_{\mathbb{R}} \mathfrak{u}= M
     \end{array}\right. \]
  with $F (Z) = \frac{a (e^{- a Z} - 1)}{1 + e^{- a Z}}$ admits a unique
  solution denoted $\mathfrak{u}_M$, which satisfies for $| M |$ small enough
  (depending on $a$)
  \[ \partial_x (F +\mathfrak{u}_M) < 0. \]
\end{lemma}

\begin{proof}
  We look at the equation
  \begin{equation}
    \left\{\begin{array}{l}
      - \partial_x \mathfrak{u}+ F\mathfrak{u}+ \frac{\mathfrak{u}^2}{2} = 0\\
      \mathfrak{u} (0) \in \mathbb{R},
    \end{array}\right. \label{profum}
  \end{equation}
  where $F = \frac{a (e^{- a Z} - 1)}{1 + e^{- a Z}}$ was introduced in the
  proof of Lemma \ref{cigaretteandkerosine} (recall that $a$ depends on time).
  A classical computation shows that
  \begin{equation}
    \mathfrak{u} (x) = \frac{\exp \left( \int_0^x F
    \right)}{\frac{1}{\mathfrak{u} (0)} - \frac{1}{2} \int_0^x \exp \left(
    \int_0^y F \right) d y} . \label{nothing}
  \end{equation}
  (with $\mathfrak{u} (x) = 0$ if $\mathfrak{u} (0) = 0$). If $\mathfrak{u}
  (0) \neq 0$, then $\mathfrak{u}$ has the same sign as $\mathfrak{u} (0)$.
  Denoting $\mathcal{F}= \exp \left( \int_0^x F \right)$, remark that
  $\mathcal{F} (0) = 1$, $\mathcal{F}$ is positive, even because $F$ is odd,
  and
  \[ \int_{\mathbb{R}} \mathcal{F}< + \infty . \]
  In particular,
  \begin{eqnarray*}
    &  & \int_{\mathbb{R}} \mathfrak{u}\\
    & = & - 2 \int_{\mathbb{R}} \partial_x \left( \ln \left(
    \frac{1}{\mathfrak{u} (0)} - \frac{1}{2} \int_0^x \mathcal{F} (y) d y
    \right) \right)\\
    & = & 2 \ln \left( \frac{1 + \frac{\mathfrak{u} (0)}{4} \int_{\mathbb{R}}
    \mathcal{F}}{1 - \frac{\mathfrak{u} (0)}{4} \int_{\mathbb{R}} \mathcal{F}}
    \right) .
  \end{eqnarray*}
  This means that $\mathfrak{u} (0) \rightarrow \int_{\mathbb{R}}
  \mathfrak{u}$ is a bijection from $\left] \frac{- 4}{\int_{\mathbb{R}}
  \mathcal{F}}, \frac{4}{\int_{\mathbb{R}} \mathcal{F}} \right[$ to
  $\mathbb{R}$. Given $M \in \mathbb{R}$, we then define $\mathfrak{u}_M$ to
  be the solution of (\ref{profum}) with $\int_{\mathbb{R}} \mathfrak{u}_M =
  M$.
  
  \
  
  We have $F' + \frac{a^2}{2} = \frac{1}{2} F^2$ and $F (0) = 0$, hence by
  (\ref{profum}) we check that $F +\mathfrak{u}_M$ satisfies the equation
  \[ \partial_x (F +\mathfrak{u}_M) = \frac{1}{2} (F +\mathfrak{u}_M)^2 -
     \frac{a^2}{2}, \]
  which implies that
  \[ (F +\mathfrak{u}_M) (x) = a \tanh (c_M - a x) \]
  with $c_M$ defined by $\mathfrak{u}_M (0) = a \tanh (c_M)$ (for $| M |$
  small enough, $\mathfrak{u}_M (0)$ is small, hence $c_M$ is well defined by
  this equation). We deduce that
  \begin{equation}
    \partial_x (F +\mathfrak{u}_M) = \frac{- a^2}{\cosh^2 (c_M - a x)} < 0.
    \label{daxonl}
  \end{equation}
\end{proof}

\subsection{Decomposition of $f$ and equation on the norm}

We denote $M = \int_{\mathbb{R}} f$ and we decompose
\[ f =\mathfrak{u}_M + \partial_y g, \]
with
\[ g (y) = \int_{- \infty}^y f -\mathfrak{u}_M . \]

First, remark that $M$ is independent of time. This is because $\partial_t
\mathbbmss{h}_{\varepsilon} = \partial_t \varepsilon \partial_{\varepsilon}
\mathbbmss{h}_{\varepsilon}$ and by Lemma \ref{reboot}, $\int_{\mathbb{R}}
\partial_{\varepsilon} \mathbbmss{h}_{\varepsilon} = 0$. Furthermore,
$\mathfrak{u}_M$ depends on time but only through $a (t)$. However, since
$\int_{\mathbb{R}} \mathfrak{u}_M = M$, we deduce that $\int_{\mathbb{R}}
\partial_t \mathfrak{u}_M = 0$. We therefore write
\[ \partial_t \mathfrak{u}_M = \partial_t a \partial_a \mathfrak{u}_M =
   \partial_t a \partial_y \mathfrak{v}_M \]
where
\[ \mathfrak{v}_M = \int_{- \infty}^y \partial_a \mathfrak{u}_M . \]
We check easily, with the explicit dependency on $a$ of $\mathfrak{u}_M$, that
$\mathfrak{v}_M$ also decays exponentially fast at $\pm \infty$ with similar
bounds as $\mathfrak{u}_M$.

\

We continue. Take some $\nu_0 > 0$ small and assume that at time $T$,
\[ \| (1 + | y |)^3 f (y, T) \|_{L^{\infty} (\mathbb{R})} + \| \partial_y f
   (., T) \|_{H^1 (\mathbb{R})} + \left| \int_{\mathbb{R}} f (y, T) \right|
   \leqslant \nu_0 . \]
Then, using the results of subsection \ref{Voyager} to estimate
$\mathfrak{u}_M$, we have
\[ | g (y) | \leqslant \frac{K \nu_0}{(1 + | y |)^2}, \]
hence $\| g \|_{L^2 (\mathbb{R})} (T) \leqslant K \nu_0$. Furthermore,
$\partial_y g = f -\mathfrak{u}_M$ and $\partial_y^2 g = \partial_y f -
\partial_y \mathfrak{u}_M$, hence
\[ \| g \|_{H^2 (\mathbb{R})} (T) \leqslant K \nu_0 . \]
That is, taking $\nu_0$ small enough, we can make $g$ as small as we want in
$H^2 (\mathbb{R})$ at the initial time $T$.

\

Now, replacing it in (\ref{will}) and integrating the equation between $-
\infty$ and $y$ implies that
\begin{eqnarray}
  &  & t^{\frac{2 \alpha}{1 + \alpha}} \partial_t g - \partial_y^2 g -
  \frac{\alpha \varepsilon}{1 + \alpha} y \partial_y g + \left(
  H_{\varepsilon} - \frac{z_c}{1 + \alpha} +\mathfrak{u}_M \right) \partial_y
  g \nonumber\\
  & + & \frac{(\partial_y g)^2}{2} - t \partial_t z_c \partial_y g
  \nonumber\\
  & + & \tilde{S} - \frac{\alpha \varepsilon}{1 + \alpha} y\mathfrak{u}_M - t
  \partial_t z_c \mathfrak{u}_M + t^{\frac{2 \alpha}{1 + \alpha}} \partial_t
  a\mathfrak{v}_M \nonumber\\
  & = & 0.  \label{angymath}
\end{eqnarray}
We recall that $\varepsilon (t) = t^{- \frac{1 - \alpha}{1 + \alpha}}$. We now
take a weight $W \in C^2 (\mathbb{R}, \mathbb{R}^{+ \ast})$ that we will
precise later on. Taking the scalar product of the equation with $g W$ leads
to
\begin{eqnarray}
  &  & t^{\frac{2 \alpha}{1 + \alpha}} \partial_t (\| g \|_{L^2 (W)}^2) + 2
  \| \partial_y g \|_{L^2 (W)}^2 \nonumber\\
  & + & \int_{\mathbb{R}} g^2 W \left( \frac{\alpha \varepsilon}{(1 +
  \alpha)} - \partial_y (H_{\varepsilon} +\mathfrak{u}_M) - \frac{\partial_y
  W}{W} \left( H_{\varepsilon} - \frac{z_c + \alpha \varepsilon y}{1 + \alpha}
  +\mathfrak{u}_M \right) \right) \nonumber\\
  & + & \int_{\mathbb{R}} g^2 W \left( - t^{\frac{2 \alpha}{1 + \alpha}}
  \frac{\partial_t W}{W} - \frac{\partial_y^2 W}{W} + t \partial_t z_c
  \frac{\partial_y W}{W} \right) \nonumber\\
  & + & \int_{\mathbb{R}} g W \left( (\partial_y g)^2 + 2 \left( \tilde{S} -
  \frac{\alpha \varepsilon}{1 + \alpha} y\mathfrak{u}_M - t \partial_t z_c
  \mathfrak{u}_M + t^{\frac{2 \alpha}{1 + \alpha}} \partial_t a\mathfrak{v}_M
  \right) \right) \nonumber\\
  & = & 0.  \label{omwtofinish}
\end{eqnarray}
Differentiating equation (\ref{angymath}) leads to
\begin{eqnarray*}
  &  & t^{\frac{2 \alpha}{1 + \alpha}} \partial_t \partial_y g - \partial_y^3
  g - \frac{\alpha \varepsilon}{1 + \alpha} \partial_y g - \frac{\alpha
  \varepsilon}{1 + \alpha} y \partial_y^2 g + \left( H_{\varepsilon} -
  \frac{z_c}{1 + \alpha} +\mathfrak{u}_M \right) \partial^2_y g\\
  & + & \partial_y (H_{\varepsilon} +\mathfrak{u}_M) \partial_y g +
  \partial_y \left( \frac{(\partial_y g)^2}{2} \right) + \partial_y \left(
  \tilde{S} - \frac{\alpha \varepsilon}{1 + \alpha} y\mathfrak{u}_M - t
  \partial_t z_c \mathfrak{u}_M + t^{\frac{2 \alpha}{1 + \alpha}} \partial_t
  a\mathfrak{v}_M \right)\\
  & = & 0.
\end{eqnarray*}
Its scalar product with $\partial_y g W$ gives us the equality
\begin{eqnarray}
  &  & t^{\frac{2 \alpha}{1 + \alpha}} \partial_t (\| \partial_y g \|_{L^2
  (W)}^2) + 2 \| \partial^2_y g \|_{L^2 (W)}^2 \nonumber\\
  & + & \int_{\mathbb{R}} (\partial_y g)^2 W \left( \frac{- \alpha
  \varepsilon}{(1 + \alpha)} + \partial_y (H_{\varepsilon} +\mathfrak{u}_M) -
  \frac{\partial_y W}{W} \left( H_{\varepsilon} - \frac{z_c + \alpha
  \varepsilon y}{1 + \alpha} +\mathfrak{u}_M \right) \right) \nonumber\\
  & + & \int_{\mathbb{R}} (\partial_y g)^2 W \left( - t^{\frac{2 \alpha}{1 +
  \alpha}} \frac{\partial_t W}{W} - \frac{\partial_y^2 W}{W} + t \partial_t
  z_c \frac{\partial_y W}{W} \right) + \int_{\mathbb{R}} W \partial_y
  ((\partial_y g)^2) \partial_y g \nonumber\\
  & + & \int_{\mathbb{R}} 2 \partial_y g W \partial_y \left( \tilde{S} -
  \frac{\alpha \varepsilon}{1 + \alpha} y\mathfrak{u}_M - t \partial_t z_c
  \mathfrak{u}_M + t^{\frac{2 \alpha}{1 + \alpha}} \partial_t a\mathfrak{v}_M
  \right) .  \label{omwtofinish2}
\end{eqnarray}
We compute, supposing that $W$ is constant outside of a compact set, that
\[ \int_{\mathbb{R}} W \partial_y ((\partial_y g)^2) \partial_y g = \frac{-
   1}{3} \int_{\mathbb{R}} \partial_y W (\partial_y g)^3 . \]
Summing (\ref{omwtofinish}) and $\lambda$ times (\ref{omwtofinish2}) for some
$\lambda > 0$ to be determined later reads
\begin{eqnarray}
  &  & t^{\frac{2 \alpha}{1 + \alpha}} \partial_t (\| g \|_{L^2 (W)}^2 +
  \lambda \| \partial_y g \|_{L^2 (W)}^2) \nonumber\\
  & + & \int_{\mathbb{R}} g^2 W D_1 + \lambda \int_{\mathbb{R}} (\partial_y
  g)^2 W D_2 \nonumber\\
  & - & \| \partial_y g \|^2_{L^2 (W)} \| g \|_{L^{\infty}} \nonumber\\
  & + & 2 \int_{\mathbb{R}} g W \left( \tilde{S} - \frac{\alpha
  \varepsilon}{1 + \alpha} y\mathfrak{u}_M - t \partial_t z_c \mathfrak{u}_M +
  t^{\frac{2 \alpha}{1 + \alpha}} \partial_t a\mathfrak{v}_M \right)
  \nonumber\\
  & + & 2 \lambda \int_{\mathbb{R}} \partial_y g W \partial_y \left(
  \tilde{S} - \frac{\alpha \varepsilon}{1 + \alpha} y\mathfrak{u}_M - t
  \partial_t z_c \mathfrak{u}_M + t^{\frac{2 \alpha}{1 + \alpha}} \partial_t
  a\mathfrak{v}_M \right) \nonumber\\
  & + & 2 \lambda \| \partial^2_y g \|_{L^2 (W)}^2 - \frac{1}{3}
  \int_{\mathbb{R}} \partial_y W (\partial_y g)^3 \nonumber\\
  & \leqslant & 0,  \label{gtfo}
\end{eqnarray}
with
\begin{eqnarray*}
  D_1 & \assign & \frac{\alpha \varepsilon}{(1 + \alpha)} - \partial_y
  (H_{\varepsilon} +\mathfrak{u}_M) - \frac{\partial_y W}{W} \left(
  H_{\varepsilon} - \frac{z_c + \alpha \varepsilon y}{1 + \alpha}
  +\mathfrak{u}_M \right)\\
  & - & t^{\frac{2 \alpha}{1 + \alpha}} \frac{\partial_t W}{W} -
  \frac{\partial_y^2 W}{W} + t \partial_t z_c \frac{\partial_y W}{W}
\end{eqnarray*}
and
\begin{eqnarray*}
  D_2 & \assign & \frac{2}{\lambda} + \frac{- \alpha \varepsilon}{(1 +
  \alpha)} + \partial_y (H_{\varepsilon} +\mathfrak{u}_M) - \frac{\partial_y
  W}{W} \left( H_{\varepsilon} - \frac{z_c + \alpha \varepsilon y}{1 + \alpha}
  +\mathfrak{u}_M \right)\\
  & - & t^{\frac{2 \alpha}{1 + \alpha}} \frac{\partial_t W}{W} -
  \frac{\partial_y^2 W}{W} + t \partial_t z_c \frac{\partial_y W}{W} .
\end{eqnarray*}
\subsection{Estimates of $D_1$ and $D_2$ and choice of $W$}\label{killkitten}

In this subsection, we choose the values of $\lambda$ and $W$ so that $D_1$
and $D_2$ are strictly positive, and satisfy some good estimates.

\subsubsection{Estimates for $y > 0$}

The goal of this subsection is to show that for $y > 0$ we have
\[ D_1 (y) \geqslant \frac{\varepsilon \alpha}{1 + \alpha} + C e^{-
   \frac{a}{2} | y |} \text{ and } D_2 (y) \geqslant 1 \]
for some constant $C > 0$ independent of $\varepsilon$. We recall that for $y
> 0$,
\[ H_{\varepsilon} (y) = h_{\varepsilon} (z_c + \varepsilon y) = h_0 (z_c +
   \varepsilon y) + F (y) + a + \varepsilon G (y), \]
hence
\begin{eqnarray*}
  D_1 & = & \varepsilon \left( \frac{\alpha}{(1 + \alpha)} - h_0' (z_c +
  \varepsilon y) \right) - \partial_y (F +\mathfrak{u}_M) - \varepsilon
  \partial_y G\\
  & - & \frac{\partial_y W}{W} \left( h_0 (z_c + \varepsilon y) -
  \frac{z_c}{1 + \alpha} + F (y) + a - \frac{\alpha \varepsilon y}{1 + \alpha}
  +\mathfrak{u}_M \right)\\
  & - & t^{\frac{2 \alpha}{1 + \alpha}} \frac{\partial_t W}{W} -
  \frac{\partial_y^2 W}{W} + t \partial_t z_c \frac{\partial_y W}{W} .
\end{eqnarray*}
In this region, we chose $W = 1$. Then,
\[ D_1 = \varepsilon \left( \frac{\alpha}{(1 + \alpha)} - h_0' (z_c +
   \varepsilon y) \right) - \partial_y (F +\mathfrak{u}_M) - \varepsilon
   \partial_y G. \]
We recall that $- h_0' (z_c + \varepsilon y) \geqslant 0$, and from
(\ref{daxonl}) we have
\[ - \partial_y (F +\mathfrak{u}_M) = \frac{a^2}{\cosh^2 (c_M - a x)} . \]
Finally, from (\ref{44}) we have
\[ \varepsilon | \partial_y G | \leqslant C_0 \left( \varepsilon^2 \ln
   \frac{1}{\varepsilon} + \varepsilon e^{- \frac{a}{2} Z} \right), \]
which implies that, for $t \geqslant T$ with $T > 0$ large enough and $| M |$
small enough (so that $c_M$ is close to $0$), we have $- \partial_y (F
+\mathfrak{u}_M) \geqslant C e^{- \frac{a}{2} | y |}$ for some $C > 0$ and
thus
\[ D_1 (y) \geqslant \frac{\varepsilon \alpha}{1 + \alpha} + C e^{-
   \frac{a}{2} | y |} . \]
Now, we have
\[ D_2 (y) = \frac{2}{\lambda} + \frac{- \alpha \varepsilon}{(1 + \alpha)} +
   \partial_y (H_{\varepsilon} +\mathfrak{u}_M), \]
and since $| \partial_y (H_{\varepsilon} +\mathfrak{u}_M) | \leqslant a^2$,
taking $\lambda$ small enough (depending only on $a$) and $t \geqslant T$ with
$T$ large enough leads to
\[ D_2 (y) \geqslant 1. \]

\subsubsection{Estimates for $y < 0$}

The goal of this subsection is to show that for $y < 0$ and fixing a well
chosen weight $W$ we have
\[ D_1 (y) \geqslant \frac{\alpha \varepsilon}{4 (1 + \alpha)} \text{ and } D_2
   (y) \geqslant 1. \]
For $y < 0$, we recall that
\[ H_{\varepsilon} (y) = h_{\varepsilon} (z_c + \varepsilon y) = h_0 (z_c +
   \varepsilon y) + F (y) - a + \varepsilon G (y), \]
hence
\begin{eqnarray*}
  D_1 & = & \varepsilon \left( \frac{\alpha}{(1 + \alpha)} - h_0' (z_c +
  \varepsilon y) \right) - \partial_y (F +\mathfrak{u}_M) - \varepsilon
  \partial_y G\\
  & - & \frac{\partial_y W}{W} \left( h_0 (z_c + \varepsilon y) -
  \frac{z_c}{1 + \alpha} + F (y) - a - \frac{\alpha \varepsilon y}{1 + \alpha}
  +\mathfrak{u}_M \right)\\
  & - & t^{\frac{2 \alpha}{1 + \alpha}} \frac{\partial_t W}{W} -
  \frac{\partial_y^2 W}{W} + t \partial_t z_c \frac{\partial_y W}{W} .
\end{eqnarray*}
Let us estimate the coefficient in factor of $\frac{\partial_y W}{W}$ in the
second line. For $y < 0$, we have $- \frac{\alpha \varepsilon y}{1 + \alpha}
\geqslant 0$, $h_0 (z_c + \varepsilon y) - \frac{z_c}{1 + \alpha} > C_0$ a
universal constant and
\[ | F (y) - a +\mathfrak{u}_M | \leqslant K e^{- \frac{a}{2} | y |} \]
if $| M |$ is small enough. Therefore, there exists $y_0$ independent of time
such that
\[ h_0 (z_c + \varepsilon y) - \frac{z_c}{1 + \alpha} + F (y) - a -
   \frac{\alpha \varepsilon y}{1 + \alpha} +\mathfrak{u}_M \geqslant
   \frac{C_0}{2} \]
for $y \leqslant y_0$, and we choose $y_0$ to be the largest value in
$\mathbb{R}^-$ such that this holds. We then define $W \in C^2$ by $W (y) = 1$
for $y \in [y_0, 0]$ and on $] - \infty, y_0]$,
\[ \frac{\partial_y W}{W} = - w (y) \varepsilon, \]
where $w$ is a $C^2$ function that satisfies $w (y) = C > 0$ if $y \in [-
\gamma \varepsilon^{- 1}, y_0 - 1]$ for some $C, \gamma > 0$ that will be
determined later on, and $w (y) = 0$ if $y \leqslant - (\gamma + 1)
\varepsilon^{- 1}$, with $| w | + | w' | \leqslant 2 C$ everywhere, and $| w'
| \leqslant 2 C \varepsilon$ if $y \leqslant y_0 - 1$. Remark in particular
that $W$ is then constant for $y \leqslant - (\gamma + 1) \varepsilon^{- 1}$,
and that this constant is uniform in time (but depends on $C$ and $\gamma$).
Indeed, for $y \leqslant y_0$ we have
\[ W (y) = \exp \left( - \varepsilon \int_{y_0}^y w \right) . \]
This also show that for $y \leqslant y_0 - 1$, we have
\[ \left| - t^{\frac{2 \alpha}{1 + \alpha}} \frac{\partial_t W}{W} -
   \frac{\partial_y^2 W}{W} + t \partial_t z_c \frac{\partial_y W}{W} \right|
   \leqslant K \varepsilon^2 . \]
We choose $\gamma$ such that for $y \leqslant - \gamma \varepsilon^{- 1}$ we
have
\[ \frac{\alpha}{(1 + \alpha)} - h_0' (z_c + \varepsilon y) \geqslant
   \frac{\alpha}{2 (1 + \alpha)} . \]
This is possible thanks to Lemma \ref{Vegas}. For $y \leqslant - \gamma
\varepsilon^{- 1}$, we have
\[ - \frac{\partial_y W}{W} \left( h_0 (z_c + \varepsilon y) - \frac{z_c}{1 +
   \alpha} + F (y) - a - \frac{\alpha \varepsilon y}{1 + \alpha}
   +\mathfrak{u}_M \right) \geqslant 0, \]
and therefore
\[ D_1 \geqslant \frac{\varepsilon \alpha}{2 (1 + \alpha)} - \partial_y (F
   +\mathfrak{u}_M) - \varepsilon \partial_y G \geqslant \frac{\alpha
   \varepsilon}{4 (1 + \alpha)} \]
if $\gamma$ is taken large enough (depending only on $\alpha, \kappa$). Now,
for $y \in [- \gamma \varepsilon^{- 1}, y_0 - 1]$, we have
\[ - \frac{\partial_y W}{W} \left( h_0 (z_c + \varepsilon y) - \frac{z_c}{1 +
   \alpha} + F (y) - a - \frac{\alpha \varepsilon y}{1 + \alpha}
   +\mathfrak{u}_M \right) \geqslant \frac{C C_0}{2} \varepsilon, \]
therefore if we take $C$ large enough we check that
\[ D_1 \geqslant \frac{\alpha \varepsilon}{4 (1 + \alpha)} \]
there as well. Now, for $y \in [y_0 - 1, 0]$, we have
\[ D_1 \geqslant - \partial_y (F +\mathfrak{u}_M) - K \varepsilon \geqslant K
   e^{- a | y |} > 0 \]
if $T$ is large enough.

\

For $D_2$, since $\left| \frac{\partial_y W}{W} \right| \leqslant 2 C
\varepsilon$, we check that taking $\lambda$ large enough and $t \geqslant T$
with $T$ large enough, we have $D_2 (y) \geqslant 1.$

\subsubsection{Summary}

With the above choices for $W$ and $\lambda$, we have
\[ D_1 \geqslant \frac{\alpha \varepsilon}{4 (1 + \alpha)} + K e^{-
   \frac{a}{2} | y |} \]
for some $K > 0$ independent of $\varepsilon$ and
\[ D_2 \geqslant 1. \]
Remark that in the case $y > 0$, we could not have choosen a similar weight
$W$, because $D_1$ contains a term $y \frac{\partial_y W}{W}$. For $y < 0,
\frac{\partial_y W}{W} < 0$ this is a positive quantity, but for $y > 0$, this
would pose an issue.

\subsection{Estimates on the source terms}\label{magitek}

We focus here on estimates on
\[ \tilde{S} - \frac{\alpha \varepsilon}{1 + \alpha} y\mathfrak{u}_M - t
   \partial_t z_c \mathfrak{u}_M + t^{\frac{2 \alpha}{1 + \alpha}} \partial_t
   a\mathfrak{v}_M \]
with $\tilde{S} = \int_{- \infty}^y \partial_t \varepsilon
\partial_{\varepsilon} \mathbbmss{h}_{\varepsilon} (z_c + \varepsilon x) d x$. We
have
\[ t \partial_t z_c = t \partial_t \varepsilon \partial_{\varepsilon} z_c =
   \frac{- 1 + \alpha}{1 + \alpha} \varepsilon \partial_{\varepsilon} z_c \]
and $\partial_t a = \frac{- 1 + \alpha}{1 + \alpha} t^{- 1} \varepsilon
\partial_{\varepsilon} a$, and we recall from (\ref{Pulverturm}) that $|
\partial_{\varepsilon} z_c | + | \partial_{\varepsilon} a | \leqslant K \left(
\ln \frac{1}{\varepsilon} \right)^2$. By the estimates on $\mathfrak{u}_M$
from Subsection \ref{Voyager}, we deduce easily that
\begin{eqnarray*}
  &  & \left| - \frac{\alpha \varepsilon}{1 + \alpha} y\mathfrak{u}_M - t
  \partial_t z_c \mathfrak{u}_M + t^{\frac{2 \alpha}{1 + \alpha}} \partial_t
  a\mathfrak{v}_M \right|\\
  & + & \left| \partial_y \left( - \frac{\alpha \varepsilon}{1 + \alpha}
  y\mathfrak{u}_M - t \partial_t z_c \mathfrak{u}_M + t^{\frac{2 \alpha}{1 +
  \alpha}} \partial_t a\mathfrak{v}_M \right) \right|\\
  & \leqslant & K \varepsilon \left( \ln \frac{1}{\varepsilon} \right)^2 e^{-
  \frac{a}{2} | y |}
\end{eqnarray*}
if $| M |$ is small enough. Concerning $\tilde{S}$, we have $| \partial_t
\varepsilon | \leqslant K t^{- 1} \varepsilon$, and using Lemma \ref{reboot},
we check that
\[ \left| \int_{- \infty}^y \varepsilon \partial_{\varepsilon}
   \mathbbmss{h}_{\varepsilon} (z_c + \varepsilon x) d x \right| \leqslant
   \frac{K \varepsilon \left( \ln \frac{1}{\varepsilon} \right)^2}{(1 + | y
   |)^{2 \alpha}} + K e^{- 2 a | y |} \]
and if $\alpha > \frac{1}{4}$ (which is needed to get enough decay in $y$),
$\tilde{S} \in L^2 (\mathbb{R})$ with
\[ \| \tilde{S} \|_{L^2 (\mathbb{R})} \leqslant \frac{K}{t} . \]
We check also, with similar arguments, that $\partial_y \tilde{S} \in L^2
(\mathbb{R})$ and
\[ \| \partial_y \tilde{S} \|_{L^2 (\mathbb{R})} \leqslant \frac{K}{t} . \]
\subsection{End of the proof of Theorem \ref{prop111}}\label{ade}

\subsubsection{Estimates on $\| g \|_{L^2 (W)}^2 + \lambda \| \partial_y g
\|_{L^2 (W)}^2$}

By Cauchy-Schwarz and the estimates on $D_1$ and $D_2$ from subsection
\ref{killkitten}, equation (\ref{gtfo}) implies for $t$ large enough the
inequality
\begin{eqnarray*}
  &  & t^{\frac{2 \alpha}{1 + \alpha}} \partial_t (\| g \|_{L^2 (W)}^2 +
  \lambda \| \partial_y g \|_{L^2 (W)}^2)\\
  & + & \int_{\mathbb{R}} g^2 W \left( \frac{\alpha \varepsilon}{4 (1 +
  \alpha)} + K e^{- 2 a | y |} \right) + (\lambda - \| g \|_{L^{\infty}}) \|
  \partial_y g \|_{L^2 (W)}^2\\
  & + & 2 \int_{\mathbb{R}} g W \left( \tilde{S} - \frac{\alpha
  \varepsilon}{1 + \alpha} y\mathfrak{u}_M - t \partial_t z_c \mathfrak{u}_M
  \right)\\
  & - & 2 \lambda \| \partial_y g \|_{L^2 (W)} \left\| \partial_y \left(
  \tilde{S} - \frac{\alpha \varepsilon}{1 + \alpha} y\mathfrak{u}_M - t
  \partial_t z_c \mathfrak{u}_M \right) \right\|_{L^2 (W)}\\
  & + & 2 \lambda \| \partial^2_y g \|_{L^2 (W)}^2 - \frac{1}{3}
  \int_{\mathbb{R}} \partial_y W (\partial_y g)^3\\
  & \leqslant & 0.
\end{eqnarray*}
By the computations of subsection \ref{magitek}, we have by Cauchy-Schwarz
that, for $\alpha > \frac{1}{4}$,
\begin{eqnarray*}
  &  & \left| \int_{\mathbb{R}} g W \left( \tilde{S} - \frac{\alpha
  \varepsilon}{1 + \alpha} y\mathfrak{u}_M - t \partial_t z_c \mathfrak{u}_M
  \right) \right|\\
  & \leqslant & K \int_{\mathbb{R}} | g | W \left( \varepsilon \left( \ln
  \frac{1}{\varepsilon} \right)^2 e^{- \frac{a}{2} | y |} + \frac{t^{- 1}}{(1
  + | y |)^{2 \alpha}} \right)\\
  & \leqslant & K \varepsilon \left( \ln \frac{1}{\varepsilon} \right)^2
  \sqrt{\int_{\mathbb{R}} g^2 W e^{- \frac{a}{2} | y |}} + t^{- 1} \| g
  \|_{L^2 (W)} .
\end{eqnarray*}
This implies that
\begin{eqnarray*}
  &  & \int_{\mathbb{R}} g^2 W \left( \frac{\alpha \varepsilon}{4 (1 +
  \alpha)} + K e^{- 2 a | y |} \right) - \left| \int_{\mathbb{R}} g W \left(
  \tilde{S} - \frac{\alpha \varepsilon}{1 + \alpha} y\mathfrak{u}_M - t
  \partial_t z_c \mathfrak{u}_M \right) \right|\\
  & \geqslant & \frac{\alpha \varepsilon}{8 (1 + \alpha)} \| g \|_{L^2 (W)}^2
  - K t^{- 2} \varepsilon^{- 1}
\end{eqnarray*}
it $\varepsilon (t)$ is small enough. We also check that
\[ \left\| \partial_y \left( \tilde{S} - \frac{\alpha \varepsilon}{1 + \alpha}
   y\mathfrak{u}_M - t \partial_t z_c \mathfrak{u}_M \right) \right\|_{L^2
   (W)} \leqslant K \varepsilon \left( \ln \frac{1}{\varepsilon} \right)^2, \]
hence
\begin{eqnarray*}
  &  & \lambda \| \partial_y g \|_{L^2 (W)}^2 - 2 \lambda \| \partial_y g
  \|_{L^2 (W)} \left\| \partial_y \left( \tilde{S} - \frac{\alpha
  \varepsilon}{1 + \alpha} y\mathfrak{u}_M - t \partial_t z_c \mathfrak{u}_M
  \right) \right\|_{L^2 (W)}\\
  & \geqslant & \frac{\lambda}{2} \| \partial_y g \|_{L^2 (W)}^2 - K
  \varepsilon^2 \left( \ln \frac{1}{\varepsilon} \right)^4 .
\end{eqnarray*}
Furthermore, by Gargliano-Nirenberg inequality, since $W$ is bounded above and
below by constants independent of time, we have
\[ \| g \|_{L^{\infty}} \leqslant \tilde{K} \| g \|_{L^2 (\mathbb{R})}^{1 / 2}
   \| \partial_y g \|_{L^2 (\mathbb{R})}^{1 / 2} \leqslant K \| g \|_{L^2
   (W)}^{1 / 2} \| \partial_y g \|_{L^2 (W)}^{1 / 2} . \]
We continue,
\[ \left| \int_{\mathbb{R}} \partial_y W (\partial_y g)^3 \right| \leqslant 2
   C \varepsilon \int_{\mathbb{R}} W (\partial_y g)^3 \leqslant 2 C
   \varepsilon \| \partial_y g \|_{L^2 (W)}^2 \| \partial_y g \|_{L^{\infty}
   (\mathbb{R})} \]
and
\[ \| \partial_y g \|_{L^{\infty} (\mathbb{R})} \leqslant K \| \partial_y g
   \|_{L^2 (W)}^{1 / 2} \| \partial^2_y g \|_{L^2 (W)}^{1 / 2} . \]
We deduce that
\begin{eqnarray*}
  &  & 2 \lambda \| \partial^2_y g \|_{L^2 (W)}^2 - \frac{1}{3} \left|
  \int_{\mathbb{R}} \partial_y W (\partial_y g)^3 \right|\\
  & \geqslant & 2 \lambda \| \partial^2_y g \|_{L^2 (W)}^2 - C \varepsilon \|
  \partial_y g \|_{L^2 (W)}^{5 / 2} \| \partial^2_y g \|_{L^2 (W)}^{1 / 2}\\
  & \geqslant & - K \varepsilon^{4 / 3} \| \partial_y g \|_{L^2 (W)}^{10 / 3}
  .
\end{eqnarray*}
Combining these estimates leads to
\begin{eqnarray*}
  &  & t^{\frac{2 \alpha}{1 + \alpha}} \partial_t (\| g \|_{L^2 (W)}^2 +
  \lambda \| \partial_y g \|_{L^2 (W)}^2)\\
  & + & \frac{\alpha \varepsilon}{8 (1 + \alpha)} (\| g \|_{L^2 (W)}^2 +
  \lambda \| \partial_y g \|_{L^2 (W)}^2) - K \| \partial_y g \|_{L^2 (W)}^{5
  / 2} \| g \|_{L^2 (W)}^{1 / 2}\\
  & - & K \varepsilon^2 \left( \ln \frac{1}{\varepsilon} \right)^4 - K
  \varepsilon^{4 / 3} \| \partial_y g \|_{L^2 (W)}^{10 / 3}\\
  & \leqslant & 0.
\end{eqnarray*}
With $\varepsilon (t) = t^{- \frac{1 - \alpha}{1 + \alpha}}$, dividing by
$t^{\frac{2 \alpha}{1 + \alpha}}$ gives us the estimate
\begin{eqnarray*}
  &  & \partial_t (\| g \|_{L^2 (W)}^2 + \lambda \| \partial_y g \|_{L^2
  (W)}^2)\\
  & + & \frac{\alpha t^{- 1}}{8 (1 + \alpha)} (\| g \|_{L^2 (W)}^2 + \lambda
  \| \partial_y g \|_{L^2 (W)}^2)\\
  & - & K t^{- \frac{2 \alpha}{1 + \alpha}} \| \partial_y g \|_{L^2 (W)}^{5 /
  2} \| g \|_{L^2 (W)}^{1 / 2} - K t^{\frac{- 4 - 2 \alpha}{3 (1 + \alpha)}}
  \| \partial_y g \|_{L^2 (W)}^{10 / 3}\\
  & - & K t^{- \frac{2}{1 + \alpha}} (\ln t)^4\\
  & \leqslant & 0,
\end{eqnarray*}
and since $\frac{- 4 - 2 \alpha}{3 (1 + \alpha)} < - 1$ and $- \frac{2}{1 +
\alpha} < - 1$, we deduce that
\[ \partial_t \left( t^{- \frac{\alpha}{8 (1 + \alpha)}} (\| g \|_{L^2 (W)}^2
   + \lambda \| \partial_y g \|_{L^2 (W)}^2) \right) \leqslant 0. \]
Therefore, if at time $T$ large enough we take $\nu_0$ small enough (depending
on $T$), then for $t \geqslant T$,
\[ \| g \|_{L^2 (W)}^2 + \lambda \| \partial_y g \|_{L^2 (W)}^2 \leqslant K
   t^{- \frac{\alpha}{8 (1 + \alpha)}} \]
for some $K > 0$ large. We check also, with similar computations, that for
some $\mu > 0$ we have
\[ \| g \|_{L^2 (W)}^2 + \lambda \| \partial_y g \|_{L^2 (W)}^2 + \mu \|
   \partial^2_y g \|_{L^2 (W)}^2 \leqslant K t^{- \frac{\alpha}{8 (1 +
   \alpha)}} . \]
This can be done by differentiating twice equation (\ref{angymath}) and taking
its scalar product with $\partial_y^2 g W$, and adding it to (\ref{gtfo}).

\subsubsection{Returning to the original scaling}

Since $g = f -\mathfrak{u}_M$, we have with $\delta = \frac{\alpha}{8 (1 + \alpha)}$ that
\[ \| f -\mathfrak{u}_M \|_{L^2 (\mathbb{R})}^2 \leqslant K \| \partial_y g
   \|_{L^2 (W)}^2 \leqslant K t^{- \delta} \]
and
\[ \| f -\mathfrak{u}_M \|_{L^{\infty} (\mathbb{R})}^2 \leqslant K \|
   \partial_y g \|_{H^1 (W)}^2 \leqslant K t^{- \delta} . \]
We recall that with $u$ solving the viscous Burgers equation, we wrote
\[ H (y, t) = t^{\frac{\alpha}{1 + \alpha}} u \left( (z_c + \varepsilon (t) y)
   t^{\frac{1}{1 + \alpha}}, t \right) \]
and
\[ H (y, t) = H_{\varepsilon (t)} (y) +\mathfrak{u}_M (y) + (f
   -\mathfrak{u}_M) (y) . \]
Therefore,
\[ \left\| t^{\frac{\alpha}{1 + \alpha}} u (x, t) - (H_{\varepsilon (t)}
   +\mathfrak{u}_M) \left( t^{\frac{1 - \alpha}{1 + \alpha}} \left( x t^{-
   \frac{1}{1 + \alpha}} - z_c (t) \right) \right) \right\|_{L^{\infty}
   (\mathbb{R})} = o_{t \rightarrow + \infty} (1), \]
and
\[ H_{\varepsilon (t)} \left( t^{\frac{1 - \alpha}{1 + \alpha}} \left( x t^{-
   \frac{1}{1 + \alpha}} - z_c (t) \right) \right) =\mathbbmss{h}_{\varepsilon
   (t)} \left( t^{- \frac{1}{1 + \alpha}} x \right) . \]
This completes the proof of Theorem \ref{prop111}.

\subsection{Proof of Theorem \ref{thfinal}.}\label{gxd}

This section is devoted to the proof of Theorem \ref{thfinal}, that is about
the equation
\[ \partial_t u - \partial_x^2 u + \partial_x \left( \frac{u^2}{2} + J (u)
   \right) = 0. \]
Doing the same change of variable as for the proof of Theorem \ref{prop111},
the only change in equation (\ref{angymath}) is that we add the term
\[ E_J \assign t^{\frac{2 \alpha}{1 + \alpha}} J \left( t^{- \frac{\alpha}{1 +
   \alpha}} (H_{\varepsilon} +\mathfrak{u}_M + \partial_y g) \right) . \]
We recall that $| J (u) | \leqslant K | u |^3$. The scalar product of $E_J$
with $g W$ can be controlled by
\begin{eqnarray*}
  &  & \left| \int_{\mathbb{R}} E_J g W \right|\\
  & \leqslant & K t^{- \frac{\alpha}{1 + \alpha}} \int_{\mathbb{R}} |
  H_{\varepsilon} +\mathfrak{u}_M + \partial_y g |^3 | g | W\\
  & \leqslant & K t^{- \frac{\alpha}{1 + \alpha}} \| g \|_{L^{\infty}} (1 +
  \| \partial_y g \|_{L^{\infty}}) \left( \| \partial_y g \|_{L^2 (W)}^2 +
  \int_{\mathbb{R}} W (H_{\varepsilon} +\mathfrak{u}_M)^2 \right) .
\end{eqnarray*}
Since $H_{\varepsilon} (y) = h_{\varepsilon} (z_c + \varepsilon y)$, we check
that if $\alpha > \frac{1}{2}$ we have
\[ \int_{\mathbb{R}} W (H_{\varepsilon} +\mathfrak{u}_M)^2 \leqslant K
   \varepsilon (t)^{- 1}, \]
hence
\begin{eqnarray*}
  &  & \left| \int_{\mathbb{R}} E_J g W \right|\\
  & \leqslant & K t^{- \frac{\alpha}{1 + \alpha}} \| g \|_{L^{\infty}} (1 +
  \| \partial_y g \|_{L^{\infty}}) \left( \| \partial_y g \|_{L^2 (W)}^2 + K
  t^{\frac{1 - \alpha}{1 + \alpha}} \right),
\end{eqnarray*}
and to consider it as an error term to conclude as in subsection \ref{ade}, we
need $t^{- \frac{2 \alpha}{1 + \alpha}} \int_{\mathbb{R}} E_J g W$ to decay in
time strictly faster than $t^{- 1 - \delta}$ for $\delta > 0$ small provided
that $\| g \|_{H^2 (\mathbb{R})}^2 \leqslant K t^{- \delta}$. This is the case
if $\frac{1 - 4 \alpha}{1 + \alpha} < - 1$, that is $\alpha > \frac{2}{3}$.

\

We can check similarly that we can treat $\int_{\mathbb{R}} \partial_y E_J
\partial_y g W$ and $\int_{\mathbb{R}} \partial_y^2 E_J \partial_y^2 g W$
similarly. For the later one, we use the fact that we also control $\mu \|
\partial_y^3 g \|_{L^2 (W)}^2 .$

\appendix\section{Proof of Proposition \ref{heatfo}}\label{sec3}

\begin{proof}
  We recall that the solution of the heat equation is
  \[ f (x, t) = \frac{1}{\sqrt{4 \pi t}} \int_{\mathbb{R}} f_0 (y) e^{-
     \frac{(x - y)^2}{4 t}} d y. \]
  We compute
  \[ t^{\alpha / 2} f \left( \sqrt{t} z, t \right) = \frac{t^{\alpha /
     2}}{\sqrt{4 \pi}} \int_{\mathbb{R}} f_0 \left( y \sqrt{t} \right) e^{- (y
     - z)^2 / 4} d y. \]
  Take any $\beta \in \left] \frac{\alpha}{2}, \frac{1}{2} \right[$. Then
  \[ t^{\alpha / 2} \left| \int_{| y | \leqslant t^{- \beta}} f_0 \left( y
     \sqrt{t} \right) e^{- (y - z)^2 / 4} d y \right| \leqslant K
     t^{\frac{\alpha}{2} - \beta} \rightarrow 0 \]
  when $t \rightarrow + \infty$. Furthermore, for $| y | \geqslant t^{-
  \beta}$, we have $| y | \sqrt{t} \geqslant t^{\frac{1}{2} - \beta}
  \rightarrow + \infty$ when $t \rightarrow + \infty$, hence
  \[ t^{\alpha / 2} f_0 \left( y \sqrt{t} \right) \rightarrow \frac{-
     \kappa}{| y |^{\alpha}} \]
  when $t \rightarrow + \infty$. We deduce that
  \[ \left| t^{\alpha / 2} f \left( \sqrt{t} z, t \right) +
     \frac{\kappa}{\sqrt{4 \pi}} \int_{| y | \geqslant t^{- \beta}} \frac{1}{|
     y |^{\alpha}} e^{- (y - z)^2 / 4} d y \right| = o_{t \rightarrow +
     \infty} (1) \]
  and since $\alpha < 1$, $| y |^{- \alpha}$ is integrable near $0$, which
  concludes the proof.
\end{proof}

\end{document}